\documentclass[11pt,bezier]{article}
\usepackage{amsmath,amssymb,amsfonts,euscript,graphicx,hyperref}
\usepackage[all]{xy}

\newtheorem{prethm}{{\bf Theorem}}

\newtheorem{prepro}{{\bf Theorem}}

\newtheorem{preprop}{{\bf Proposition}}

\newenvironment{prop}{\begin{preprop}{\hspace{-0.5
               em}{\bf.}}}{\end{preprop}}

\newtheorem{precor}{{\bf Corollary}}

\newenvironment{cor}{\begin{precor}{\hspace{-0.5
               em}{\bf.}}}{\end{precor}}

\newtheorem{preconj}{{\bf Conjecture}}

\newtheorem{preremark}{{\bf Remark}}

\newenvironment{remark}{\begin{preremark}\rm{\hspace{-0.5
               em}{\bf.}}}{\end{preremark}}

\newtheorem{preexample}{{\bf Example}}

\newenvironment{example}{\begin{preexample}\rm{\hspace{-0.5
               em}{\bf.}}}{\end{preexample}}

\newtheorem{prelem}{{\bf Lemma}}

\newtheorem{prelam}{{\bf Lemma}}

\newtheorem{preproof}{{\bf Proof.}}

\begin{document}
\title{\bf \large Dessins on Modular Curves}
\author{{\large{Khashayar Filom}}\thanks{Department of Mathematical Sciences, Sharif University of Technology, Tehran, Iran.\hspace{3cm} E-mail: $\mathsf{filom@mehr.sharif.ir}$ The work of this author is part of his
 Ph.D. thesis at Sharif University under co-advising of H. Fanai and M. Shahshahani.} \quad {\large{Ali Kamalinejad}}\thanks{School of Mathematics, Institute for Research in Fundamental Sciences (IPM), Tehran, Iran. \hspace{3cm}E-mail: $\mathsf{kamalinejad@ipm.ir}$}}

\date{}
\maketitle
\begin{abstract}
Given a finite index subgroup $\Gamma$ of ${\rm{PSL}}_2(\Bbb{Z})$, we investigate Belyi functions on the corresponding modular curve $X(\Gamma)$ by introducing two methods for constructing such functions. Numerous examples have been worked out completely and as an application, we have derived modular equations for $\Gamma_0(2),\Gamma_0(3)$ and several special values of the $j$-function by a new method based on the theory of Belyi functions and dessin d'enfants.
\end{abstract}

\section{Introduction}
 It is an elementary fact that in the action of $\Gamma(1)={\rm{SL}}_2(\Bbb{Z})$ on the upper half plane $\Bbb{H}$, there are two elliptic orbits, that of ${\rm{i}}$ and ${\rm{e}}^{\frac{\pi{\rm{i}}}{3}}$. So for any finite index subgroup $\Gamma$ of  $\Gamma(1)$
the finite ramified cover $\Gamma\backslash\Bbb{H}\rightarrow\Gamma(1)\backslash\Bbb{H}$ after compactification gives rise to a morphism from the compact Riemann surfaces associated with $\Gamma$ to the compactification of $\Gamma(1)\backslash\Bbb{H}$ with at most three ramification values which are $\Gamma(1).{\rm{i}}$ and  $\Gamma(1).{\rm{e}}^{\frac{\pi{\rm{i}}}{3}}$ along with $\Gamma(1).\infty$. On the other hand, the  $j$-function establishes an isomorphism from the compactification of $\Gamma(1)\backslash\Bbb{H}$ to $\Bbb{CP}^1$. Therefore, the aforementioned morphism of compact Riemann surfaces can be realized as a Belyi function.\\
The purpose of this paper is to show how dessin theory can be an effective tool in computational problems in modular
curves. More specifically, it allows the explicit computation of modular equations, the determination of certain special values of $j$-function and provides explicit examples of ramified coverings of algebraic curves. 
In particular, in the case of isogenies of elliptic curves, these examples lead to identities for the Weierstrass $\wp$- functions, for instance in Remark \ref{identity}, that appear to be new. \\
Examples of calculations of the Belyi function and the dessin d'enfant associated with a congruence subgroup have appeared before in the literature e.g. in \cite{Mckay2},\cite{Mckay1} for the case of genus zero torsion-free congruence subgroups. We will investigate some of the congruence  subgroups which possess torsion elements for instance in Examples \ref{Gamma_0(2)}, \ref{Gamma_0(3)} and furthermore some subgroups of higher genera in 
Examples \ref{Gamma(6)}, \ref{Gamma(8)} and \ref{Gamma_0(11)}. Moreover, a second method for constructing Belyi functions on modular curves is introduced in section 3, the advantage of which is a simple description of the monodromy representation. The final section is devoted to applications to computing the modular equations and some special $j$-values. \\
For notations we mainly follow \cite{Milne} which is our primary reference on modular functions and modular curves. $\Bbb{H}^*:=\Bbb{H}\cup\Bbb{Q}\cup\{\infty\}$ denotes the upper half plane union cusps for $\Gamma(1)$. $\Gamma(N)$ is the principal congruence subgroup of level $N$ and $\Gamma_0(N)$ is the subgroup of $\Gamma(1)={\rm{SL}}_2(\Bbb{Z})$ consisting  of matrices $\begin{bmatrix}
a&b\\
c&d
\end{bmatrix}$
with $c\equiv 0 \pmod N$. $\Gamma$ always indicates a finite index subgroup of $\Gamma(1)$. We denote the compact Riemann surface $\Gamma\backslash\Bbb{H}^*$ by $X(\Gamma)$. For $\Gamma=\Gamma(N)$ or $\Gamma_0(N)$, notations $X(N)$ and $X_0(N)$ are used instead. The point ${\rm{e}}^{\frac{\pi{\rm{i}}}{3}}$ of the upper half plane is denoted by $\rho$ and $[z]$ means the orbit of $z\in\Bbb{H}$ under the action of $\Gamma$ and also the corresponding point in $X(\Gamma)$. Finally, $\bar{\Gamma}$ is the image of $\Gamma\leq{\rm{SL}_2(\Bbb{Z})}$ in $\bar{\Gamma}(1)={\rm{PSL}}_2(\Bbb{Z})\subset{\rm{Aut}}(\Bbb{H})$ and we will work with M\"obius transformations $T:z\mapsto z+1$ and $S:z\mapsto\frac{-1}{z}$ (induced by matrices 
$T=\begin{bmatrix}
1 & 1\\
0 &1
\end{bmatrix} $ and 
$S=\begin{bmatrix}
0 & -1\\
1 &0
\end{bmatrix} $
 respectively) as a set of generators for ${\rm{PSL}}_2(\Bbb{Z})$. 

\section{First Construction}
Let us analyze the ramification structure of the map $X(\Gamma)=\Gamma\backslash\Bbb{H}^*\rightarrow X(1)=\Gamma(1)\backslash\Bbb{H}^*$. The commutative diagram below
$$
\xymatrix{\Bbb{H} \ar[d]   \ar[dr]& \\
                 \Gamma\backslash\Bbb{H}\ar[r] & \Gamma(1)\backslash\Bbb{H}
}
$$
and the fact that in the action of $\bar{\Gamma}(1)$ on $\Bbb{H}$ points with non-trivial stabilizer are exactly the elements of orbits $\Gamma(1).\rho$ or  $\Gamma(1).{\rm{i}}$ with stabilizers of orders three and two respectively, indicates that branch values of the map $\Gamma\backslash\Bbb{H}\rightarrow\Gamma(1)\backslash\Bbb{H}$ obtained from its restriction are among $[\rho],[{\rm{i}}]\in\Bbb{H}$. Moreover, the point $[z]\in\Gamma\backslash\Bbb{H}$  where $z\in\Gamma(1).\rho$ (resp. $z\in\Gamma(1).{\rm{i}}$) is a ramification point of this map iff $z$ is not an elliptic point of $\Gamma$ and in that case, multiplicity of this point is three (resp. two). The only other ramification value that $\Gamma\backslash\Bbb{H}^*\rightarrow\Gamma(1)\backslash\Bbb{H}^*$ may possess is $[\infty]$ whose fiber is the set of orbits of cusps for $\Gamma$. Hence $\Gamma\backslash\Bbb{H}^*\rightarrow\Gamma(1)\backslash\Bbb{H}^*$  is a Belyi function provided one identifies $X(1)=\Gamma(1)\backslash\Bbb{H}^*$ with $\Bbb{CP}^1$ in a way that the subset containing two elliptic orbits (points $[\rho],[{\rm{i}}]\in X(1)$) and one orbit of cusps (point $[\infty]\in X(1)$) bijects to $\{0,1,\infty\}$. The modular function $j:\Bbb{H}\rightarrow\Bbb{C}$ gives us such an identification since it has a simple pole at infinity and satisfies $j(\rho)=0$ and $j({\rm{i}})=1728$. Thus, for any $\Gamma$, the function $\frac{1}{1728}{j}:X(\Gamma)=\Gamma\backslash\Bbb{H}^*\rightarrow\Bbb{CP}^1$ is Belyi:
\begin{prop}\label{firstBelyi}
Let $\Gamma$ be a finite index subgroup of $\Gamma(1)$. Then 
$\begin{cases}
f:X(\Gamma)=\Gamma\backslash\Bbb{H}^*\rightarrow\Bbb{CP}^1\\
[z]\mapsto\frac{1}{1728}j(z)
\end{cases}$
is a Belyi function of degree $m:=[\bar{\Gamma}(1):\bar{\Gamma}]$ whose black  (resp. white) vertices are $\Gamma$-orbits of points in $\Gamma(1).\rho$  (resp. $\Gamma(1).{\rm{i}}$) and centers of faces ($\times$ vertices) are orbits of cusps of $\Gamma$. If $\nu_3$ (resp. $\nu_2$) and $\nu_\infty$ indicate the number of inequivalent elliptic points of order three (resp. two) and the number of inequivalent cusps for $\Gamma$ respectively, then the dessin has $\nu_\infty$ faces and among its black (resp. white) vertices, there are $\frac{m-\nu_3}{3}$ (resp. $\frac{m-\nu_2}{2}$) vertices of degree three (resp. two) and the rest,  i.e. $\nu_3$ (resp. $\nu_2$) remaining black (resp. white) vertices, are all of degree one\footnote {Using these details to calculate the Euler characteristic of the triangulation obtained from the dessin gives rise to the expression $1+\frac{m}{12}-\frac{\nu_2}{4}-\frac{\nu_3}{3}-\frac{\nu_\infty}{2}$ for the genus of $X(\Gamma)$ as in \cite{Milne}.}. Furthermore, the number of edges of the face surrounding a vertex $[z]$ of type $\times$ where $z\in\Bbb{Q}\cup\{\infty\}$ is twice the width of the cusp $z$ of $\Gamma$.     
\end{prop}  
This description of degrees of vertices of various kinds in our dessin certainly does not need to specify  isomorphism class of  the corresponding Belyi function since there is still the question of the way the graph is embedded on the underlying topological surface. This deficiency is remedied in the second method for constructing Belyi functions on modular curves introduced in the next section where in certain cases the monodromy of Belyi function can be recovered from the corresponding subgroup of $\Gamma(1)$. Nevertheless, symmetries stemming from the group action always can be used to determine the dessin. In particular, when  $\Gamma\unlhd\Gamma(1)$ the well-defined  transitive left action of $\Gamma(1)$ on $\Gamma\backslash\Bbb{H}^*$  implies that the deck transformation group of the Belyi function $\Gamma\backslash\Bbb{H}^*\rightarrow\Gamma(1)\backslash\Bbb{H}^*\cong\Bbb{CP}^1$ may be identified with ${\bar{\Gamma}(1)}/{\bar{\Gamma}}$. Consequently, this ramified covering is regular, our dessin is uniform (degrees of vertices of the same type coincide) and ${\bar{\Gamma}(1)}/{\bar{\Gamma}}$ is actually a group of symmetries of the dessin. Exploiting these symmetries extremely simplifies calculations with dessins as we shall see in examples of this section where we mainly deal with principal congruence subgroups $\Gamma(N)$:
\begin{cor}\label{principalcongruence}
  Since for  $N>1$, $\Gamma(N)$ does not have any elliptic point, according to
 \ref{firstBelyi}, in their dessins, $\bullet,\circ,\times$ vertices are of degrees $2,3,N$, respectively. The genus of $X(N)$ is given by  $1+\frac{N-6}{12N}\mu_N$ where $\mu_N:=\left[\bar{\Gamma}(1):\bar{\Gamma}(N)\right]=\begin{cases} 6\qquad\qquad\qquad\qquad \quad N=2\\ \frac{N^3}{2}\prod_{\substack{p\mid N\\ p \text{ prime} }}(1-\frac{1}{p^2})\quad N>2\end{cases}$ cf. \cite{Milne}.  As a result, there are $\frac{\mu_N}{3}$ black vertices of degree three, $\frac{\mu_N}{2}$ white vertices of degree two, $\mu_N$ edges and $\frac{\mu_N}{N}$ faces each with $2N$ edges.       
\end{cor}

\begin{remark}\label{complexstructure}
Although, as mentioned before, the mere description of dessin's degree sequence in $\ref{firstBelyi}$ is not sufficient to determine the corresponding Belyi function $X(\Gamma)\rightarrow\Bbb{CP}^1$, it might be enough to determine the complex structure of the compact Riemann surface $X(\Gamma)$ provided that this degree sequence specifies the underlying graph of the dessin up to isomorphism. This is due to a simple observation: various embeddings of a bipartite graph as a dessin on a fixed oriented topological surface all lead to the same complex structure. The reason is that different embeddings differ by an orientation preserving self-homeomorphism of the surface. Such a homeomorphism respects two triangulations of the surface obtained from these embeddings and thus respects complex structures on the surface induced by these triangulations. There are some relevant cases where the underlying graph is uniquely determined, e.g. when the dessin is uniform.     
 To see this, note that since each white vertex is of degree two, the dessin is obtained from a graph with black vertices by putting one white vertex on each edge\footnote{In literature dessins of this kind are called {\it{clean}}. It follows from Proposition \ref{firstBelyi} that the dessin on $X(\Gamma)$ is clean provided that $\Gamma$ does not have any elliptic point of order two, e.g. in the case of $\Gamma=\Gamma(N)$ when $N\geq2$ or $\Gamma_0(N)$  when $ -1$  is not a quadratic residue modulo $N$.}. This is a $3$-regular graph embedded (as a map) on a topological surface where  every face has the same number of edges and the number of vertices, edges and faces are prescribed by $\Gamma$. Clearly such a graph is unique up to isomorphism. \end{remark}

\begin{example}\label{Gamma(2)}
For our first example,  consider the index six normal subgroup $\Gamma(2)$. Since $\mu_2=6$,  Corollary \ref{principalcongruence} implies that the genus of $X(2)$ is zero, there are $\mu_2=6$ edges, $\frac{\mu_2}{3}=2$ black vertices each of degree three, $\frac{\mu_2}{2}=3$ white vertices each of degree two and $\frac{\mu_2}{N}=3$ faces each with $2N=4$ edges. Thus, the dessin is an embedding of the graph $K_{2,3}$ into the sphere. There is an isomorphism  ${\bar{\Gamma}(1)}/{\bar{\Gamma}(2)}\stackrel{\cong}{\rightarrow}{\rm{S}}_3$ in which cosets of  $S$ and $T$ are mapped to transpositions. Using these symmetries, after fixing three white vertices by $ [{\rm{i}}]=0, [{\rm{i}}+1]=1,  \left[\frac{{\rm{i}}+1}{2}\right]=-1$, one can quickly recover the Belyi function $f:X(2)\rightarrow X(1)\cong\Bbb{CP}^1$ along with its dessin depicted in Figure 1:
\begin{equation} \label{Gamma(2)Belyi}
f(z)=\frac{(3z^2+1)^3}{(9z^2-1)^2}\quad f(z)-1=\frac{27z^2(z^2-1)^2}{(9z^2-1)^2} 
\end{equation}
\noindent We may also explicitly describe the group ${\rm{Deck}}\left(f:\Bbb{CP}^1\rightarrow\Bbb{CP}^1\right)\cong {\rm{S}}_3$ as the group generated by order two elements $z\mapsto -z$ and $z\mapsto \frac{-z+1}{3z+1}$ induced by the left actions of $S$ and $T$ respectively\footnote{The degree six meromorphic function on the  Riemann sphere obtained as $\Gamma(2)\backslash\Bbb{H}^*\rightarrow\Gamma(1)\backslash\Bbb{H}^*$ is in fact up to a constant factor isomorphic to the formula of the $j$-invariant of an elliptic curve relative to its Legendre form which is given by $\lambda\mapsto 256\frac{(\lambda^2-\lambda+1)^3}{(\lambda^2-\lambda)^2}$, cf. page 267 of \cite{dessin}.}.  
\begin{figure}[ht]
\centering
\includegraphics[height=4.5cm]{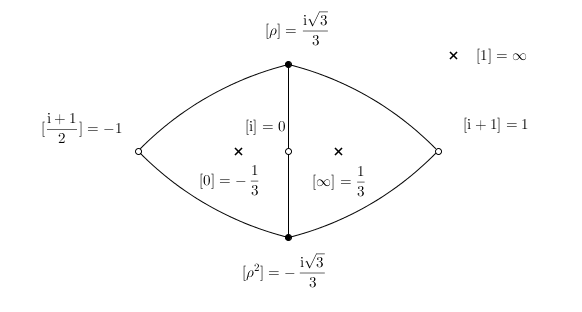}\\
Figure 1
\end{figure}
\end{example}

\begin{remark}\label{draw}
Although irrelevant in the previous example, it might seem ambiguous how one should connect black and white vertices in the dessin on $X(\Gamma)$. Note that in the identification of $X(1)$ with $\Bbb{CP}^1$ in Proposition \ref{firstBelyi} via a multiple of $j$-invariant, the arc $\left\{z\in\Bbb{H}\big|\, |z|=1,0\leq\mathfrak{R}(z)\leq\frac{1}{2}\right\}$ from $\rho$ to $\rm{i}$ in the boundary of the usual fundamental domain for $\Gamma(1)={\rm{SL}_2}(\Bbb{Z})$ bijects to the unit interval. Consequently, edges of dessin which are its preimages can be described as $\left\{\left[\gamma.{\rm{e}}^{{\rm{i}}\theta}\right] \big | \frac{\pi}{3}\leq \theta\leq\frac{\pi}{2},\gamma\in {\rm{SL}_2}(\Bbb{Z})\right\}$. For torsion-free $\Gamma$ there is a general method discussed in \cite{Mckay2} for drawing the underlying graph of the dessin by relating it to the {\it{Schreier coset  graph}} associated with the subgroup $\Gamma$ of ${\rm{PSL}}_2(\Bbb{Z})$.
\end{remark}

\begin{example}\label{Gamma(3)}
In this example, we concentrate on the principal congruence subgroup $\Gamma(3)$ and the corresponding Belyi function $f:X(3)\rightarrow X(1)\cong\Bbb{CP}^1$ where $N=3, \mu_3=12$. So by Corollary \ref{principalcongruence}, the corresponding dessin has twelve edges, four black vertices each of degree three, six white vertices each of degree two and four faces ($\times$ vertices) each surrounded with six edges.  The deck transformation group $\bar{\Gamma}(1)/\bar{\Gamma}(3)$ is isomorphic to the alternating group ${\rm{A}}_4$ by the map which sends cosets of $T$ and $S$  to a $3$-cycle and a product of two disjoint transpositions in ${\rm{A}}_4$, respectively.  The description $\left(\Bbb{Z}_2\times\Bbb{Z}_2\right)\rtimes\Bbb{Z}_3$ of the group ${\rm{A}}_4$ as a semidirect product of the subgroup of order two elements  by a subgroup generated by an order three element (like  the coset of $T$) indicates that in terms of $S$ and $T$ a set of representatives  may be described as: 
\begin{equation} \label{representatives}
\{I_2,S,TST^{-1},T^2ST^{-2},T,TS,T^2ST^{-1},ST^{-2},T^2,T^2S,ST^{-1},TST^{-2}\}
\end{equation}
Applying them to points $\rho ,{\rm{i}},\infty$ of $\Bbb{H}^*$, one can easily compute the $\Gamma(3)$-orbits that form the set of vertices  and then Remark \ref{draw} enables us to draw the underlying graph of our dessin on $X(3)\cong\Bbb{CP}^1$ as in Figure 2.  
\begin{figure}[ht]
\centering
\includegraphics[height=4.5cm]{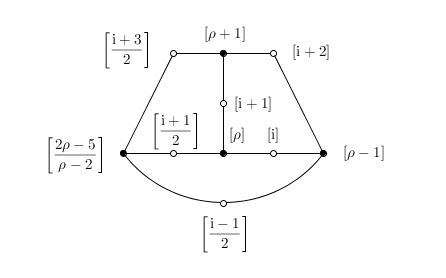}\\
Figure 2
\end{figure}

\noindent $T$ is mapped to a $3$-cycle in ${\rm{A}}_4$ and hence acts as an order three M\"obius transformation of $\Bbb{CP}^1$ whose action on either of black vertices or $\times$ vertices has exactly one fixed point. But $T$ leaves $\Gamma(3)$-orbits of $\left[\frac{2\rho-5}{\rho-2}\right]$ and $[\infty]$ of   invariant. So killing three degrees of freedom in drawing a graph as a dessin on the Riemann sphere, we may rigidify the graph in Figure 2 by assuming that $T$ is the $120\,^\circ$ rotation about the origin whose fixed points are $0$ and $\infty$, i.e. $T$ induces the deck transformation $z\mapsto {\rm{e}}^{\frac{2\pi{\rm{i}}}{3}}z$ and furthermore $\left[\frac{2\rho-5}{\rho-2}\right]=0,[\infty]=\infty$. 
After setting $\alpha=[{\rm{i}}]\,\,\beta=\left[\frac{{\rm{i}}+3}{2}\right]$, coordinates of remaining white vertices will be: 
\begin{equation*}
\begin{split}
&\left[{\rm{i}}+1\right]=T\left(\left[{\rm{i}}\right]\right)={\rm{e}}^{\frac{2\pi{\rm{i}}}{3}}\alpha\quad \left[{\rm{i}}+2\right]=T^2\left(\left[{\rm{i}}\right]\right)={\rm{e}}^{\frac{4\pi{\rm{i}}}{3}}\alpha\\ &\left[\frac{{\rm{i}}-1}{2}\right]=T\left(\left[\frac{{\rm{i}}+3}{2}\right]\right)={\rm{e}}^{\frac{2\pi{\rm{i}}}{3}}\beta\quad \left[\frac{{\rm{i}}+1}{2}\right]=T^2\left(\left[\frac{{\rm{i}}+3}{2}\right]\right)={\rm{e}}^{\frac{4\pi{\rm{i}}}{3}}\beta
\end{split}
\end{equation*}
On the other hand:
 $$S\left(\left[{\rm{i}}+1\right] \right)=\left[S({\rm{i}}+1)=\frac{{\rm{i}}-1}{2}\right]\quad S\left(\left[\frac{{\rm{i}}+1}{2}\right]\right)=\left[S\left(\frac{{\rm{i}}+1}{2}\right)={\rm{i}}-1\right]=\left[{\rm{i}}+2\right]$$ 
which implies that the involution of $\Bbb{CP}^1$ induced by $S$ interchanges ${\rm{e}}^{\frac{2\pi{\rm{i}}}{3}}\alpha,{\rm{e}}^{\frac{4\pi{\rm{i}}}{3}}\alpha$ and also ${\rm{e}}^{\frac{2\pi{\rm{i}}}{3}}\beta,{\rm{e}}^{\frac{4\pi{\rm{i}}}{3}}\beta$ while fixing points $\alpha$ and $\beta$, the observation that yields its formula as $z\mapsto\frac{(\alpha+\beta)z-2\alpha\beta}{2z-(\alpha+\beta)}$ and requires $\alpha,\beta\in\Bbb{C}-\{0\}$ to satisfy $\alpha^2+4\alpha\beta+\beta^2=0$. Next task is to determine the coordinates of vertices of type $\bullet$ or $\times$ in terms of $\alpha$ and $\beta$.  In its action on black vertices, $T$ remains $\left[\frac{2\rho-5}{\rho-2}\right]=0$ invariant while $S$, that acts by the M\"obius transformation $z\mapsto\frac{(\alpha+\beta)z-2\alpha\beta}{2z-(\alpha+\beta)}$ on $X(3)\cong\Bbb{CP}^1$, takes it to $[\rho+1]$  which gives us $\frac{2\alpha\beta}{\alpha+\beta}=-(\alpha+\beta)$ for the coordinate of the black vertex $[\rho+1]$. Now since $T$, which induces $z\mapsto {\rm{e}}^{\frac{2\pi{\rm{i}}}{3}}z$, permutes other three black vertices  as $[\rho-1]\mapsto[\rho]\mapsto[\rho+1]\mapsto[\rho+2]=[\rho-1]$, we deduce that $[\rho-1]=-{\rm{e}}^{\frac{2\pi{\rm{i}}}{3}}(\alpha+\beta), [\rho]=-{\rm{e}}^{\frac{4\pi{\rm{i}}}{3}}(\alpha+\beta)$. With the same argument, for vertices of type $\times$, $[\infty]=\infty$ yields: $[0]=S\left([\infty]\right)=\frac{(\alpha+\beta)z-2\alpha\beta}{2z-(\alpha+\beta)}\big |_{\infty}=\frac{\alpha+\beta}{2}$ and from that: $[1]=T\left([0]\right)={\rm{e}}^{\frac{2\pi{\rm{i}}}{3}}\frac{\alpha+\beta}{2},[2]=T^2\left([0]\right)={\rm{e}}^{\frac{4\pi{\rm{i}}}{3}}\frac{\alpha+\beta}{2}$. Summarizing all these calculations, we arrive at coordinates of all vertices in terms of $\alpha$ and $\beta$:
\begin{equation} \label{list3}
\begin{split}
&\bullet\,\text{vertices}: \left[\frac{2\rho-5}{\rho-2}\right]=0, [\rho+1]=-(\alpha+\beta),[\rho-1]=-{\rm{e}}^{\frac{2\pi{\rm{i}}}{3}}(\alpha+\beta), [\rho]=-{\rm{e}}^{\frac{4\pi{\rm{i}}}{3}}(\alpha+\beta)\\
& \circ\,\text{vertices}: [{\rm{i}}]=\alpha, \left[{\rm{i}}+1\right]={\rm{e}}^{\frac{2\pi{\rm{i}}}{3}}\alpha, \left[{\rm{i}}+2\right]={\rm{e}}^{\frac{4\pi{\rm{i}}}{3}}\alpha, \left[\frac{{\rm{i}}+3}{2}\right]=\beta, \left[\frac{{\rm{i}}-1}{2}\right]={\rm{e}}^{\frac{2\pi{\rm{i}}}{3}}\beta, \left[\frac{{\rm{i}}+1}{2}\right]={\rm{e}}^{\frac{4\pi{\rm{i}}}{3}}\beta \\
&\times\,\text{vertices}: [\infty]=\infty, [0]=\frac{\alpha+\beta}{2}, [1]={\rm{e}}^{\frac{2\pi{\rm{i}}}{3}}\frac{\alpha+\beta}{2}, [2]={\rm{e}}^{\frac{4\pi{\rm{i}}}{3}}\frac{\alpha+\beta}{2}
\end{split}
\end{equation}
The Belyi function $f$ on the Riemann sphere whose dessin has coordinates written in \eqref{list3}, must satisfy the identities:
\begin{equation} \label{formula1}
f(z)=k\frac{z^3\left(z^3+(\alpha+\beta)^3\right)^3}{\left(z^3-\left(\frac{\alpha+\beta}{2}\right)^3\right)^3}\quad f(z)-1=k\frac{\left(z^3-\alpha^3\right)^2\left(z^3-\beta^3\right)^2}{\left(z^3-\left(\frac{\alpha+\beta}{2}\right)^3\right)^3} 
\end{equation}
After some rather complicated calculations one verifies that aside from $\alpha^2+4\alpha\beta+\beta^2=0$, the only other constraint imposed by above identities is $k=\frac{1}{8(\alpha+\beta)^3}$, yielding $f(z)=\frac{1}{8(\alpha+\beta)^3}\frac{z^3\left(z^3+(\alpha+\beta)^3\right)^3}{\left(z^3-\left(\frac{\alpha+\beta}{2}\right)^3\right)^3}$. In  other words, we have the following interesting identity whenever $\alpha^2+4\alpha\beta+\beta^2=0$:
$$ \frac{1}{8\left(\alpha+\beta\right)^3}t\left(t+(\alpha+\beta)^3\right)^3-\left(t-\left(\frac{\alpha+\beta}{2}\right)^3\right)^3= \frac{1}{8\left(\alpha+\beta\right)^3}\left(t-\alpha^3\right)^2\left(t-\beta^3\right)^2$$
Without any loss of generality, we may work with specific values $\alpha=1+\sqrt{3},\beta=1-\sqrt{3}$. Plugging them into \eqref{list3},\eqref{formula1} and the formula derived for the action of $S$ completes this example: the Belyi function $f:\Bbb{CP}^1\rightarrow\Bbb{CP}^1$ satisfies 
\begin{equation}\label{Gamma(3)Belyi}
f(z)=\frac{1}{64}\frac{z^3(z^3+8)^3}{(z^3-1)^3}\quad f(z)-1=\frac{1}{64}\frac{(z^6-20z^3-8)^2}{(z^3-1)^3}
\end{equation}
whose dessin is illustrated in Figure 3 and its group of deck transformations is generated by $z\mapsto{\rm{e}}^{\frac{2\pi{\rm{i}}}{3}}z$ coming form the left action of $T$ and $z\mapsto\frac{z+2}{z-1}$ coming from the left action of $S$ which under some suitable isomorphism ${\rm{Deck}}\left(f:\Bbb{CP}^1\rightarrow\Bbb{CP}^1\right)\stackrel{\cong}{\rightarrow}{\rm{A}}_4$ are mapped to a $3$-cycle and a product of two disjoint transpositions, respectively. 
\begin{figure}[ht]
\centering
\includegraphics[height=8.5cm]{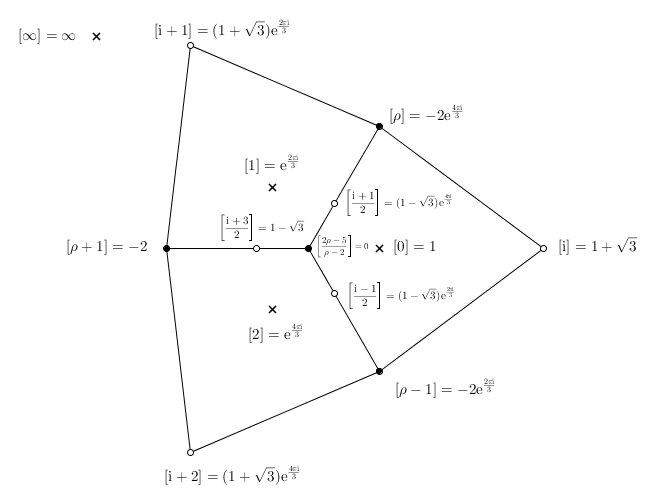}\\
Figure 3
\end{figure} 
\end{example}

\begin{example}\label{Gamma(6)}
Let us consider $N=6$ and examine the Belyi function  $f:X(6)\rightarrow X(1)\cong\Bbb{CP}^1$ which  is of degree $\mu_6=72$, corresponding to a dessin on the torus with $72$ edges, $24$ black vertices of degree three, $36$ white vertices of degree  two and $12$ faces each with $2N=12$ edges. It factorizes as 
$$X(6)=\Gamma(6)\backslash\Bbb{H}^*\stackrel{6\text{-sheeted}}{\longrightarrow}X(3)=\Gamma(3)\backslash\Bbb{H}^*\cong\Bbb{CP}^1\stackrel{12\text{-sheeted}}{\longrightarrow}\Gamma(1)\backslash\Bbb{H}^*\cong\Bbb{CP}^1$$
where the latter map is the Belyi function studied in Example \ref{Gamma(3)}. We will concentrate on the former  which expresses a genus one Riemann surface as a 6-fold ramified cover of $\Bbb{CP}^1$.  Ramification values of $\tilde{f}:X(6)=\Gamma(6)\backslash\Bbb{H}^*\rightarrow X(3)=\Gamma(3)\backslash\Bbb{H}^*$ are the points corresponding to cusps of $\Gamma(3)$, that  are  $[0],[1],[2],[\infty]\in X(3)$ and the fiber above each has precisely three points all of multiplicity two. The  deck transformation group
 is a version of ${\rm{S}}_3$ generated by order two and three automorphisms  of $X(6)=\Gamma(6)\backslash\Bbb{H}^*$ induced respectively by  left actions of $T^3$ and 
$\begin{bmatrix}
1 & 3\\
3 & 10
\end{bmatrix}$. We will employ them to construct factorizations of $\tilde{f}$ that not only will lead to its equation, but also provide an interesting explicit example of an isogeny between elliptic curves. First, note that any intermediate covering  in $\tilde{f}:X(6)\rightarrow X(3)$ satisfies the same property as $\tilde{f}$, that is multiplicities of its ramification points are all equal to two, the feature that simplifies applying the Riemann-Hurwitz formula: in any factorization  $X(6)\stackrel{2\text{-sheeted}}{\longrightarrow}Y\stackrel{3\text{-sheeted}}{\longrightarrow} X(3)$ of the 6-fold map $\tilde{f}:X(6)\rightarrow X(3)$, the genus of $Y$ is zero whereas in any factorization in the form of $X(6)\stackrel{3\text{-sheeted}}{\longrightarrow}Z\stackrel{2\text{-sheeted}}{\longrightarrow} X(3)$, the genus of $Z$ is one. With the help of order two and three deck transformations induced by
$T^3=
\begin{bmatrix}
1 & 3\\
0 & 1 
\end{bmatrix}$
 and  
$\begin{bmatrix}
1 & 3\\
3 & 10
\end{bmatrix}$,
factorizations of $\tilde{f}$ of both kinds just introduced  can be exhibited:
\begin{equation*}
\begin{split}
&\Gamma(6)\backslash\Bbb{H}^*\stackrel{2-\text{sheeted}}{\longrightarrow}\left(\Gamma(6).\left\langle \begin{bmatrix}
1 & 3\\
0 & 1 
\end{bmatrix}\right\rangle\right)\backslash\Bbb{H}^*\stackrel{3-\text{sheeted}}{\longrightarrow}\Gamma(3)\backslash\Bbb{H}^*\\
&\Gamma(6)\backslash\Bbb{H}^*\stackrel{3-\text{sheeted}}{\longrightarrow}\left(\Gamma(6).\left\langle \begin{bmatrix}
1 & 3\\
3 & 10
\end{bmatrix}\right\rangle\right)\backslash\Bbb{H}^*\stackrel{2-\text{sheeted}}{\longrightarrow}\Gamma(3)\backslash\Bbb{H}^*
\end{split}
\end{equation*}
They fit into the commutative diagram below:
\begin{equation} \label{diagram1}
\xymatrixcolsep{7pc}\xymatrix{ 
X(6)=\Gamma(6)\backslash\Bbb{H}^*\ar[d]_{2-\text{sheeted}}  \ar[r]^{3-\text{sheeted}} \ar[dr]_{\tilde{f} \text{ (6-sheeted)}} & \left(\Gamma(6).\left\langle {\begin{bmatrix}
1 & 3\\
3 & 10
\end{bmatrix}}\right\rangle\right)\backslash\Bbb{H}^* \ar[d]^{2-\text{sheeted}}\\
\left(\Gamma(6).\left\langle {\begin{bmatrix}
1 & 3\\
0 & 1 
\end{bmatrix}}\right\rangle\right)\backslash\Bbb{H}^*  \ar[r]_{3-\text{sheeted}}^{x\mapsto h(x)}      & X(3)=\Gamma(3)\backslash\Bbb{H}^*}
\end{equation}
 where genera of the surfaces in the top and bottom rows are one and zero respectively. Thus, in columns we have ramified 2-fold covers that allow us to write $\Gamma(6)\backslash\Bbb{H}^*$ and  $\left(\Gamma(6).\left\langle \begin{bmatrix}
1 & 3\\
3 & 10
\end{bmatrix}\right\rangle\right)\backslash\Bbb{H}^*$ as elliptic curves in the hyperelliptic form $y^2=(x-\alpha)(x-\beta)(x-\gamma)$ and each column as the projection $(x,y)\mapsto x$. The top row will be a degree three isogeny between these elliptic curves. By diagram chasing one deduces that the degree three map $h:\Bbb{CP}^1\rightarrow\Bbb{CP}^1$ has exactly four branch values which form the set of the branch values of both $\tilde{f}$ and  right column of \eqref{diagram1}. Moreover, the fiber of $h$  over any of them consists of a point of multiplicity two and a simple point. These simple points of $h$ over its branch values  are  the branch values of the left column. Identifying the modular curve $X(3)$ with the Riemann sphere just like Example \ref{Gamma(3)}, the coordinate of cusps or in other words, branch values of $\tilde{f}:X(6)\rightarrow X(3)$ are: $ [\infty]=\infty, [0]=1, [1]={\rm{e}}^{\frac{2\pi{\rm{i}}}{3}}, [2]={\rm{e}}^{\frac{4\pi{\rm{i}}}{3}}$ cf. Figure 3. Hence  the curve in the top right corner of the diagram \eqref{diagram1} is $y^2=x^3-1$ and we will find the other elliptic curve via exhibiting a degree three meromorphic function $h$ on the Riemann sphere with the aforementioned ramification structure. The function $h(x)=\frac{4-x^3}{3x^2}$ has the desired properties. It attains the value $1$ at points $-2,1$ where the latter is simple and due to the way that $h$ transforms under $120\,^\circ$ rotation, similarly  ${\rm{e}}^{\frac{2\pi{\rm{i}}}{3}},{\rm{e}}^{\frac{4\pi{\rm{i}}}{3}}$ are its non-critical points above critical values ${\rm{e}}^{\frac{2\pi{\rm{i}}}{3}},{\rm{e}}^{\frac{4\pi{\rm{i}}}{3}}$. So the critical values of the left column are also $0,1,{\rm{e}}^{\frac{2\pi{\rm{i}}}{3}},{\rm{e}}^{\frac{4\pi{\rm{i}}}{3}}$ and hence the other elliptic curve  may  be described by $y^2=x^3-1$ as well! Now, in the the top row of \eqref{diagram1} the isogeny may be described as $(x,y)\mapsto \left(h\left(x\right),g\left(x\right)y\right)$ where $g(x)$ must satisfy $g(x)^2=\frac{h(x)^3-1}{x^3-1}$ which yields: $g(x)=\pm\frac{{\rm{i}}\sqrt{3}}{9}\frac{x^3+8}{x^3}$. Thus, our computations culminate at an explicit description of a degree three self-isogeny of the elliptic curve $y^2=x^3-1$:
$\begin{cases}
\left\{y^2=x^3-1\right\}\rightarrow \left\{y^2=x^3-1\right\}\\
(x,y)\mapsto \left(\frac{4-x^3}{3x^2},\frac{{\rm{i}}\sqrt{3}}{9}\frac{x^3+8}{x^3}y\right)
\end{cases}$. Combining this isogeny with the 2-fold  ramified covering $(x,y)\mapsto x$ yields the equation of $\tilde{f}$  as $(x,y)\mapsto\frac{4-x^3}{3x^2}$. Composing it with $X(3)\cong\Bbb{CP}^1\rightarrow X(1)\cong\Bbb{CP}^1$ derived in  \eqref{Gamma(3)Belyi} gives rise to the following equation for the degree $72$ Belyi  function $f:X(6)\rightarrow X(1)\cong\Bbb{CP}^1$ and concludes this example:\footnote{Here, the Belyi function on the elliptic curve factors through a ramified 2-fold cover from the elliptic curve to $\Bbb{CP}^1$. Such a property of a dessin on a hyperelliptic Riemann surface is called {\it{properness}} and was first introduced in \cite{Kamalinejad}. It is a powerful tool in computations with dessin d'enfants, cf. \cite{cwd}.}
\begin{equation} \label{Gamma(6)Belyi} 
\begin{cases}
\left\{y^2=x^3-1\right\}\rightarrow\Bbb{CP}^1\\
(x,y)\mapsto -\left[\frac{(x^3-4)\left((x^3-4)^3-216x^6\right)}{12x^2(x^3-1)(x^3+8)^2}\right]^3
\end{cases}
\footnote{Therefore $x\mapsto -\left[\frac{(x^3-4)\left((x^3-4)^3-216x^6\right)}{12x^2(x^3-1)(x^3+8)^2}\right]^3$ is the Belyi function  associated with the congruence subgroup 
$\Gamma(6).\left\langle
\begin{bmatrix}
1 & 3\\
0& 1
\end{bmatrix}
 \right\rangle$. So according to the table on page 275 of \cite{Mckay2} of conjugacy classes of  genus zero, torsion-free congruence subgroups of index 36 and their ramification data, this group coincides with $\Gamma_0(2)\cap\Gamma(3)$.}
\end{equation}
\end{example}

\begin{remark}\label{identity}
The explicit equation $(x,y)\mapsto\frac{4-x^3}{3x^2}$ for a degree three self-isogeny of the hexagonal elliptic curve $\{y^2=x^3-1\}$ results in the following interesting identity involving the Weierstrass function $\wp$ associated with the lattice $\Bbb{Z}+\Bbb{Z}\omega$ where $\omega={\rm{e}}^{\frac{2\pi{\rm{i}}}{3}}$:
$$
\wp\left(\left(1+2\omega\right)z\right)-\wp\left(\frac{1}{2}\right)+\frac{\wp\left(\frac{\omega}{2}\right)-\wp\left(\frac{1}{2}\right)}{\omega-1}=\frac{4\left(\frac{\wp\left(\frac{\omega}{2}\right)-\wp\left(\frac{1}{2}\right)}{\omega-1}\right)^3-\left(\wp(z)-\wp\left(\frac{1}{2}\right)+\frac{\wp\left(\frac{\omega}{2}\right)-\wp\left(\frac{1}{2}\right)}{\omega-1}\right)^3}{3\left(\wp(z)-\wp\left(\frac{1}{2}\right)+\frac{\wp\left(\frac{\omega}{2}\right)-\wp\left(\frac{1}{2}\right)}{\omega-1}\right)^2}
$$
 This idea can be applied to any self-isogeny in order to derive non-trivial identities for Weierstrass functions. For instance, there is a degree five self-isogeny  of the square elliptic curve\footnote{The existence of these isogenies means that there is a trisection of the hexagonal tessellation of plane by the same tessellation and also the square tessellation of plane may be decomposed to five versions of itself. Realizing these observations geometrically is a challenging exercise. }:
$$\begin{cases}
\left\{y^2=x^3-x\right\}\rightarrow \left\{y^2=x^3-x\right\}\\
(x,y)\mapsto \left(\frac{(-1+2{\rm{i}})^2x(x^2-(1+2{\rm{i}}))^2}{(5x^2-(1-2{\rm{i}}))^2},\frac{(-1+2{\rm{i}})^3(x^2-(1+2{\rm{i}}))(x^4+(8{\rm{i}}+2)x^2+1)}{(5x^2-(1-2{\rm{i}}))^3}y\right)
\end{cases}$$
which yields the following identity for the Weierstrass function of the lattice $\Bbb{Z}+\Bbb{Z}{\rm{i}}$:
\begin{equation*}
\wp\left(\left(1+2{\rm{i}}\right)z\right)=(-1+2{\rm{i}})^2\frac{\wp(z)\left(\wp(z)^2-(1+2{\rm{i}})\wp\left(\frac{{\rm{i}}}{2}\right)^2\right)^2}{\left(5\wp(z)^2-(1-2{\rm{i}})\wp\left(\frac{{\rm{i}}}{2}\right)^2\right)^2}
\end{equation*}
\end{remark}
\section{Second Construction}
In our second construction, we assume $\Gamma$ to be a finite index subgroup of $\Gamma(2)$ instead of $\Gamma(1)={\rm{PSL}}_2(\Bbb{Z})$. It is no loss of generality as the Belyi theorem implies that arithmetic Riemann surfaces are exactly those  containing a proper Zariski open subset uniformised by a finite index subgroup of $\Gamma(2)$. This is due to a key observation: $\Bbb{H}\rightarrow\Gamma(2)\backslash\Bbb{H}$ is a universal covering since $\Gamma(2)$ does not possess any elliptic element and on the other hand $\Gamma(2)\backslash\Bbb{H}$ may be identified with $\Bbb{C}-\{0,1\}$ because it is the complement of cusps' orbits in $X(2)=\Gamma(2)\backslash\Bbb{H}^*\cong\Bbb{CP}^1$  while the  number of inequivalent  cusps for $\Gamma(2)$ is exactly three:
\begin{equation}\label{cusps}
[0]=\left\{\frac{a}{b}\mid a\, \text{even}, b\, \text{odd} \right\}\quad [1]=\left\{\frac{a}{b}\mid a\, \text{odd}, b\, \text{odd} \right\}\quad [\infty]=\left\{\frac{a}{b}\mid a\, \text{odd}, b\, \text{even} \right\}\cup\{\infty\}
\end{equation}
So in this section, we fix an isomorphism from $X(2)$ to $\Bbb{CP}^1$  which takes points $[0],[1],[\infty]\in\Gamma(2)\backslash\Bbb{H}^*$ to $0,1,\infty\in\Bbb{CP}^1$, respectively.  Then for any finite index subgroup $\Gamma$ of $\Gamma(2)$, the obvious map  $X(\Gamma)=\Gamma\backslash\Bbb{H}^*\rightarrow X(2)=\Gamma(2)\backslash\Bbb{H}^*\cong\Bbb{CP}^1$ of degree $\left[\bar{\Gamma}(2):\bar{\Gamma}\right]$ may be assumed to be Belyi. In the corresponding dessin, vertices are orbits of cusps of $\Gamma$ where the degree of a vertex is half the width of the corresponding cusp:
\begin{equation}\label{vertices2}
\begin{cases}
\bullet\,\text{vertices}\\
\circ\,\text{vertices}\\
\times\,\text{vertices}
\end{cases}:
\text{ orbits that partition the } \Gamma  \text{-invariant set}
\begin{cases}
\left\{\frac{a}{b}\mid a\, \text{even}, b\, \text{odd} \right\}\\
\left\{\frac{a}{b}\mid a\, \text{odd}, b\, \text{odd} \right\}\\
\left\{\frac{a}{b}\mid a\, \text{odd}, b\, \text{even} \right\}\cup\{\infty\}
\end{cases}  
\end{equation}
The main advantage of working with this Belyi function rather than $X(\Gamma)\rightarrow X(1)\cong\Bbb{CP}^1$ introduced in the previous section, aside from lower degree that eases computations, is that the monodromy homomorphism of $X(\Gamma)\rightarrow X(2)\cong\Bbb{CP}^1$, which specifies the isomorphism class of the Belyi function, may be investigated via knowledge of the subgroup $\Gamma$ as we shall explain.
Away from the set $\{[0],[1],[\infty]\}$ containing its ramification values, $X(\Gamma)\rightarrow X(2)$ restricts to $\Gamma\backslash\Bbb{H}\rightarrow\Gamma(2)\backslash\Bbb{H}$ which is an intermediate cover of the universal cover $\Bbb{H}\rightarrow\Gamma(2)\backslash\Bbb{H}$\footnote{This is the difference between this approach and that of the second section where the covering $\Bbb{H}\rightarrow\Gamma(1)\backslash\Bbb{H}$ is ramified.}. Thus, $\bar{\Gamma}(2)=\Gamma(2)/\{\pm{\rm{I}}_2\}$, which acts freely on $\Bbb{H}$, can be identified with the  rank two free group $\pi_1\big(\Gamma(2)\backslash\Bbb{H}\cong\Bbb{C}-\{0,1\}\big)$  and this intermediate cover corresponds to the subgroup $\bar{\Gamma}$ of $\bar{\Gamma}(2)$. We fix the following two generators for the free group $\bar{\Gamma}(2)=\Gamma(2)/\{\pm I_2\}$:
$$T^2=\begin{bmatrix}
1&2\\
0&1
\end{bmatrix}\quad ST^2S^{-1}=-ST^2S=
\begin{bmatrix}
1&0\\
-2&1
\end{bmatrix}$$
The monodromy is a homomorphism from this rank two free group  to the permutations of a regular fiber. In general, given a Belyi function, fixing a basis for the free group $\pi_1\big(\Bbb{C}-\{0,1\}\big)$, it is a standard fact that the monodromy homomorphism is determined with a pair of permutations of a regular fiber that generate a transitive subgroup and such pairs modulo simultaneous conjugation  are  in a one-to-one correspondence with isomorphism classes of  Belyi functions. See \cite{dessin} for more details and also a very concrete description of the monodromy of a Belyi function in terms of its dessin. Here, as the left action $\bar{\Gamma}(2)$ on $\Bbb{H}$ is free, fibers may be identified with the set of right cosets of $\bar{\Gamma}$ in $\bar{\Gamma}(2)$: points above $[x]\in\Gamma(2)\backslash\Bbb{H}$ are precisely $[\beta_1.x],\dots,[\beta_m.x]\in\Gamma\backslash\Bbb{H}$ where $\{\beta_1,\dots,\beta_m\}$  is a set of right coset representatives of $\bar{\Gamma}$ in $\bar{\Gamma}(2)$. Let us analyze the monodromy homomorphism. A fundamental domain for $\Gamma(2)$ is depicted in Figure 4 as union of six fundamental domains for $\Gamma(1)=\rm{SL}_2(\Bbb{Z})$. 
\begin{figure}[ht]
\centering
\includegraphics[height=6cm,width=12cm]{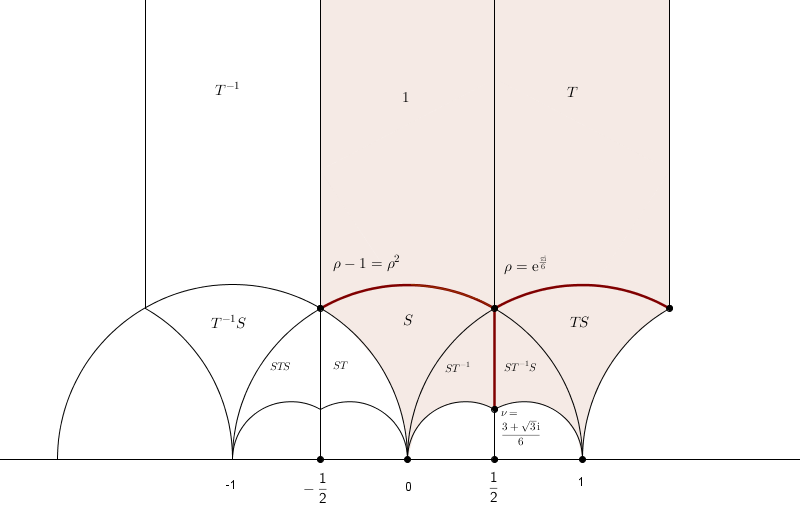}\\
Figure 4
\end{figure}
\noindent Two arcs and the segment which are bolded in the figure provide us with two small positive loops in $\Gamma(2)\backslash\Bbb{H}$ based at $[\rho-1]$ around $[\infty]$ and $[0]$ whose homotopy classes form a basis for the free group $\pi_1\big(\Gamma(2)\backslash\Bbb{H},[\rho-1]\big)$: by applying the quotient map to the path traversing first the bolded circular arc from $\rho-1$ to $\rho$ and then the bolded circular arc from $\rho$ to $\rho+1$ we will get a small counterclockwise loop in the punctured Riemann sphere $\Gamma(2)\backslash\Bbb{H}$ based at $[\rho-1]=[\rho+1]$ around the cusp $[\infty]$. Similarly, if one traverses the bolded segment from  $\nu=\frac{3+\sqrt{3}{\rm{i}}}{6}$ to $\rho$ and then the bolded circular arc from $\rho$ to $\rho-1$  and maps the resulting path to $\Gamma(2)\backslash\Bbb{H}$ via $\Bbb{H}\rightarrow\Gamma(2)\backslash\Bbb{H}$, a small counterclockwise loop  in $\Gamma(2)\backslash\Bbb{H}$ based at $[\nu]=[\rho-1]$ around $[0]$ will be determined. Fixing the basis for $\pi_1\big(\Gamma(2)\backslash\Bbb{H},[\rho-1]\big)$ consisting of these two homotopy classes, it is easy to describe the monodromy homomorphism from  $\pi_1\big(\Gamma(2)\backslash\Bbb{H},[\rho-1]\big)$ to the permutation group on elements of the fiber of $\Gamma\backslash\Bbb{H}\rightarrow\Gamma(2)\backslash\Bbb{H}$  above $[\rho-1]$. Since the endpoints of the liftings  to the universal cover $\Bbb{H}$ just constructed    satisfy $\rho+1=T^2(\rho-1)$ and $\rho-1=ST^2S^{-1}(\nu)$, we conclude that the corresponding two permutations on the fiber above $[\rho-1]$, identified with the set of right cosets of $\bar{\Gamma}$ in $\bar{\Gamma}(2)$, are those induced by right multiplication with $T^2$ and $ST^2S^{-1}$ \footnote{Equivalently, having in mind the way we identified $X(2)$ with $\Bbb{CP}^1$ and following notations of \cite{dessin}, in the permutation representation $(\sigma_0,\sigma_1,\sigma_\infty)$ (where $\sigma_0\sigma_1\sigma_\infty=1$) of the Belyi function $X(\Gamma)\rightarrow X(2)\cong\Bbb{CP}^1$, $\sigma_0$ is the permutation induced by multiplication with $ST^2S^{-1}$ and $\sigma_\infty$ is the permutation induced by multiplication with $T^2$.}.
This discussion proves:
\begin{prop}\label{secondBelyi}
Let $\Gamma$ be a finite index subgroup of $\Gamma(2)$. Then the obvious map 
$$f: X(\Gamma)=\Gamma\backslash\Bbb{H}^*\rightarrow X(2)=\Gamma(2)\backslash\Bbb{H}^*\cong\Bbb{CP}^1$$
where $X(2)$ is identified with $\Bbb{CP}^1$ in the way that $[0]\mapsto 0$, $[1]\mapsto 1$ and $[\infty]\mapsto\infty$, 
is a Belyi function of degree $m:=[\bar{\Gamma}(2):\bar{\Gamma}]$. In its dessin black, white and vertices  corresponding to faces are $\Gamma$-orbits of points in $\Gamma(2).0$, $\Gamma(2).1$ and $\Gamma(2).\infty$ respectively. The number of  edges is $m$ and the degree of a vertex $[z]\in X(\Gamma)$ is half the width of the cusp $z\in\Bbb{Q}\cup\{\infty\}$ of $\Gamma$. Moreover, the monodromy is given by two permutations of the set of right cosets of $\bar{\Gamma}$ in $\bar{\Gamma}(2)$ induced by the right actions of $T^2$  and $ST^2S^{-1}$.
\end{prop}
\noindent When $\Gamma\unlhd\Gamma(2)$, the Belyi function is regular and its dessin is uniform and by \ref{vertices2} the number of $\bullet,\circ$ or $\times$ vertices  are respectively the number of $\Gamma$-orbits that subsets $\Gamma(2).0$, $\Gamma(2).1$ or $\Gamma(2).\infty$ of $\Bbb{Q}^*=\Bbb{Q}\cup\{\infty\}$ decompose to. Using  normality, this number is easy to find: it only suffices to consider the index in $\bar{\Gamma}(2)$ of the subgroup generated by $\bar{\Gamma}$ and stabilizers of $0$, $1$ or $\infty$ in the action of $\bar{\Gamma}(2)$:   
 \begin{cor}\label{secondBelyisymmetric}
Suppose $\Gamma$ is a finite index normal subgroup of $\Gamma(2)$. Then the Belyi function
$f: X(\Gamma)\rightarrow X(2)$ from Proposition \ref{secondBelyi} is regular.
Put:
\begin{equation*}
\begin{split}
&m=\left[\bar{\Gamma}(2):\bar{\Gamma}\right]\quad k_0={\rm{min}}\left\{ l\in\Bbb{N}\mid \left(ST^2S^{-1}\right)^l\in\pm\Gamma\right\}\\
& k_1={\rm{min}}\left\{ l\in\Bbb{N}\mid \left((TS)T^2(TS)^{-1}\right)^l\in\pm\Gamma\right\}\quad k_\infty={\rm{min}}\left\{ l\in\Bbb{N}\mid T^{2l}\in\pm\Gamma\right\}
\end{split}
\end{equation*}
Then the degrees of $\bullet$, $\circ$ or $\times$ vertices, in the same order, are given by $k_0$, $k_1$ and $k_\infty$ and the numbers of vertices of these type are $\frac{m}{k_0}$, $\frac{m}{k_1}$ and $\frac{m}{k_\infty}$. Hence, the genus of $X(\Gamma)$  equals: $$g=1+\frac{m}{2}\left(1-\frac{1}{k_0}-\frac{1}{k_1}-\frac{1}{k_\infty}\right)$$ \\
The monodromy representation of this degree $m$ Belyi function $f$ is given by two permutations of ${\rm{S}}_m$ induced by the right actions of $T^2$  and $ST^2S^{-1}$ on the underlying set of the order $m$ group $\bar{\Gamma}(2)/\bar{\Gamma}$. \\   
\end{cor}  
Considering the case of principal congruence subgroups $\Gamma(2N)$ which are normal in $\Gamma(1)$ while contained in $\Gamma(2)$, the preceding corollary indicates that: 
\begin{cor}\label{secondprincipalcongruence}
The genus of $X(2N)$ equals $1+\frac{m}{2}\left(1-\frac{3}{N}\right)$ where:
$$m=\left[\bar{\Gamma}(2):\bar{\Gamma}(2N)\right]=\begin{cases}
4\quad \quad \! N=2\\
\mu_N\quad N\geq 3 \text{ odd}\\
\frac{4}{3}\mu_N\, \,N\geq 4\text{ even}
\end{cases} $$   
For the dessin on $X(2N)$ constructed in Proposition \ref{secondBelyi},  the number of edges of the dessin is $m$, each black or white vertex is of degree $N$ and each face is surrounded with $2N$ edges. The monodromy  is determined by two permutations of the set of elements of the order $m$ group 
$ \bar{\Gamma}(2)/\bar{\Gamma}(2N)\cong {\rm{ker}}\left({\rm{PSL}}_2(\Bbb{Z}_{2N})\rightarrow{\rm{PSL}}_2(\Bbb{Z}_{2})\right)$
induced by right  multiplication with elements $T^2=\begin{bmatrix}
1&2\\
0&1
\end{bmatrix}, ST^2S^{-1}=
\begin{bmatrix}
1&0\\
-2&1
\end{bmatrix}$ of $\Gamma(1)={\rm{SL}}_2(\Bbb{Z})$.
\end{cor}

\begin{example}\label{Gamma(4)}
Let $N=2$. According to Corollary \ref{secondprincipalcongruence} the Belyi function $f:X(4)\cong\Bbb{CP}^1\rightarrow X(2)\cong\Bbb{CP}^1$  is of degree four and in  its dessin  every vertex is of degree two.  The quotient group $ \bar{\Gamma}(2)/\bar{\Gamma}(4)$ is isomorphic to
to the Klein four-group. Corollary \ref{secondprincipalcongruence} implies that the monodromy homomorphism is specified by  a pair of permutations on the set of elements of Klein group  induced by  multiplication with two non-trivial elements of this group. So after an arbitrary labeling of four edges, monodromy  is given by two distinct elements of the symmetric group ${\rm{S}}_4$ each a product of two disjoint transpositions, e.g. $(1\,2)(3\,4)$ and $(1\,4)(2\,3)$. The  the dessin is illustrated completely in the Figure 5 which (after arbitrarily fixing three vertices) immediately yields the equation of the Belyi function $f$:
\begin{figure}[ht]
\centering
\includegraphics[height=4.5cm]{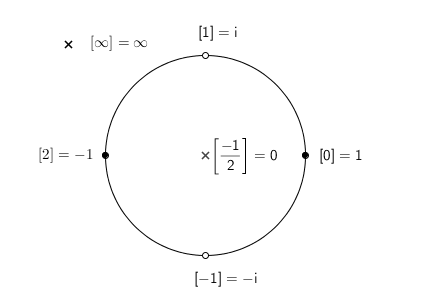}\\
Figure 5
\end{figure} 
\begin{equation}\label{Gamma(4)Belyi}
f(z)=-\frac{1}{4}\frac{(z^2-1)^2}{z^2}\quad f(z)-1=-\frac{1}{4}\frac{(z^2+1)^2}{z^2}
\end{equation}
 The deck transformation group of the Belyi function $f$  is generated by holomorphic involutions $z\mapsto -z$  and $z\mapsto\frac{1}{z}$  corresponding to  left actions of $T^2$ and $ST^2S^{-1}$ on $\Bbb{H}^*$, respectively.  Finally, with substituting $z$ in the expression \eqref{Gamma(4)Belyi} of this Belyi function   with $\frac{z+1}{3(z-1)}$, in order to change the identification of $X(2)$ with $\Bbb{CP}^1$ introduced in the beginning of this section to the one used in Example \ref{Gamma(2)} and Figure 1, and then composing it with the Belyi function $X(2)\rightarrow X(1)\cong\Bbb{CP}^1$ calculated in \ref{Gamma(2)Belyi}, one derives the equation below for the degree $24$ Belyi function $X(4)\rightarrow X(1)\cong\Bbb{CP}^1$ (and its difference with one) suggested in our first construction:
$$\frac{1}{108}\frac{(z^8+14z^4+1)^3}{z^4(z^4-1)^4}-1=\frac{1}{108}\frac{(z^4+1)^2(z^4+6z^2+1)^2(z^4-6z^2+1)^2}{z^4(z^4-1)^4}$$
\end{example}

\begin{example}\label{Gamma(8)}
Here we sketch a similar analysis for the case of $N=4$ in Corollary \ref{secondprincipalcongruence} where we have the Belyi function $f:X(8)\rightarrow X(2)\cong\Bbb{CP}^1$  of degree $\frac{4}{3}\mu_4=32$ on the Riemann surface $X(8)$ of genus five. This Belyi function  has exactly $\frac{32}{4}=8$ ramification points all of multiplicity $N=4$ in each of its critical fibers. But this map factors through the degree four Belyi function $X(4)\rightarrow X(2)$ from the previous example whose ramification points are of multiplicity two (note that all of these coverings are regular). Hence, it  suffices only to analyze the ramified $8$-fold cover $\tilde{f}:X(8)\rightarrow X(4)\cong\Bbb{CP}^1$ whose ramification points are  of multiplicity two and its branch values are the ramification points of 
$\begin{cases}X(4)\rightarrow X(2)\\ 
z\mapsto -\frac{1}{4}\frac{(z^2-1)^2}{z^2} 
\end{cases}$ that appeared in \eqref{Gamma(4)Belyi}. We conclude that the branch values of $\tilde{f}$ are $\pm 1,\pm{\rm{i}},0,\infty$ and the fiber above each of them consists of four ramification points of multiplicity two. The group of deck transformations of $\tilde{f}$ is the elementary abelian $2$-group $ \bar{\Gamma}(4)/\bar{\Gamma}(8)\cong{\rm{ker}}\left({\rm{PSL}}_2(\Bbb{Z}_{8})\rightarrow{\rm{PSL}}_2(\Bbb{Z}_{4})\right)\cong\Bbb{Z}_2\times\Bbb{Z}_2\times\Bbb{Z}_2$.
The cosets of the matrices $T^4=\begin{bmatrix} 1&4\\0&1 \end{bmatrix}$, $ST^4S^{-1}=\begin{bmatrix} 1&0\\-4&1\end{bmatrix}$
and $A:=(TS)T^4(TS)^{-1}=\begin{bmatrix} -3&4\\-4&5 \end{bmatrix}$ in $\bar{\Gamma}(4)/\bar{\Gamma}(8)$ form a basis for this group over $\Bbb{Z}_2$. Now the order two deck transformation of $\tilde{f}$ induced by the left action of an element $\gamma\in\Gamma(4)-\Gamma(8)$,  in its action on a  critical fiber (which consists of four critical points) either fixes all the points or permutes them as a product of two disjoint transpositions. The former happens only if $\gamma$ is congruent to a parabolic matrix modulo eight.  This is the case for the elements of the cosets of $T^4,ST^4S^{-1},A$ while the action of  any other $\gamma$ on critical fibers of $\tilde{f}$ is fixed point free and thus results in the  factorization $X(8)\rightarrow \langle\gamma\rangle\backslash X(8)\rightarrow X(4)$ of $\tilde{f}:X(8)\rightarrow X(4)$  where the first map is an unramified $2$-fold covering a genus three Riemann surface by the genus five Riemann surface $X(8)$. Next, let us consider the deck transformations induced by $T^4,ST^4S^{-1},A$. They fix some of the ramification points of $\tilde{f}$ because they are parabolic matrices. They are conjugate in ${\rm{SL}}_2(\Bbb{Z})$ and therefore by analyzing just one of them it can be easily verified that each of these three deck transformations fixes all the points in exactly two of critical fibers and in any of other four critical fibers permutes four points of that fiber as a product of two disjoint transpositions. Using the Riemann-Hurwitz formula, we deduce that for any $\gamma$ in $\{T^4,ST^4S^{-1},A\}.\Gamma(8)$, $X(8)\rightarrow \langle\gamma\rangle\backslash X(8)$ is a degree two ramified covering of a genus one Riemann surface by $X(8)$ with precisely eight ramification points all of multiplicity two. We will utilize these observation to investigate $\tilde{f}:X(8)=\Gamma(8)\backslash\Bbb{H}^*\rightarrow X(4)=\Gamma(4)\backslash\Bbb{H}^*$ via factorizing it by introducing some suitable intermediate subgroups between $\Gamma(8)$ and $\Gamma(4)$ as in the following diagram:
\begin{equation} \label{diagram3}
\xymatrixcolsep{2.7pc}\xymatrix{ 
X(8)=\Gamma(8)\backslash\Bbb{H}^*\ar[d]_{\text{unramified}}  \ar[r]  & \left(\Gamma(8).\langle A\rangle\right)\backslash\Bbb{H}^* \ar[d] & \\
\left(\Gamma(8).\langle AST^4S^{-1}\rangle\right)\backslash\Bbb{H}^*  \ar[r]      & \left(\Gamma(8).\langle A,ST^4S^{-1}\rangle\right)\backslash\Bbb{H}^* \ar[r]&  X(4)=\Gamma(4)\backslash\Bbb{H}^*\cong\Bbb{CP}^1
}\end{equation}
In this diagram, all of the arrows are of degree two and the multiplicity of any of their ramification points is two since $\tilde{f}:X(8)\rightarrow X(4)$ satisfies the same property. Only the left column is unramified, as indicated in the diagram. Hence since $\left(\Gamma(8).\langle A\rangle\right)\backslash\Bbb{H}^*$  is of genus one as discussed before, we conclude that $\left(\Gamma(8).\langle A,ST^4S^{-1}\rangle\right)\backslash\Bbb{H}^*$ is of genus zero and consequently the right column and the bottom row realize the genus one Riemann surface $\left(\Gamma(8).\langle A\rangle\right)\backslash\Bbb{H}^*$ and the genus three Riemann surface $\left(\Gamma(8).\langle AST^4S^{-1}\rangle\right)\backslash\Bbb{H}^*$ as hyperelliptic curves. The $2$-fold ramified cover $\left(\Gamma(8).\langle A,ST^4S^{-1}\rangle\right)\backslash\Bbb{H}^*\cong\Bbb{CP}^1\rightarrow X(4)=\left((\Gamma(8).\langle A,ST^4S^{-1},T^4\rangle=\Gamma(4)\right)\backslash\Bbb{H}\cong\Bbb{CP}^1$ that appeared on the right hand side of \eqref{diagram3} may be assumed to be $z\mapsto z^2$ because its branch values are exactly those branch values of $\tilde{f}$ where the action of the parabolic element $T^4$ on the corresponding fiber is trivial and these are points $[\infty]=\infty,[\frac{-1}{2}]=0$ of $X(4)$ (cf. Figure 5). In this situation the branch values of the right column and the bottom row are in the preimage of the set $\{\pm 1,\pm{\rm{i}},0,\infty\}$ under this map but cannot be branch values of $z\mapsto z^2$. Finding these branch values enables us to write the equations of these hyperelliptic curves explicitly and simplify the diagram \eqref{diagram3} as:
\begin{equation*}
\xymatrixcolsep{7pc}\xymatrix{X(8)\ar[d]_{\text{unramified} \, 2\text{-fold cover}} \ar[r]^{\text{ramified} \, 2\text{-fold cover}}_{\text{with eight critical points}} & \{y^2=z^4-1\}\ar[d]^{(y,z)\mapsto z} & \\ \{x^2=z^8-1\}\ar[r]_{(x,z)\mapsto z} & \Bbb{CP}^1\ar[r]_{z\mapsto z^2} & X(4)\cong\Bbb{CP}^1}
\end{equation*}
Consequently, the genus five Riemann surface $X(8)$ is birational to the affine algebraic curve $\left\{(x,y,z)\in\Bbb{C}^3\mid x^2=z^8-1,y^2=z^4-1\right\}$ on which  the degree eight ramified cover $\tilde{f}:X(8)\rightarrow X(4)\cong\Bbb{CP}^1$  restricts to $(x,y,z)\mapsto z^2$. The group ${\rm{Deck}}\left(\tilde{f}:X(8)\rightarrow X(4)\right)$ isomorphic to $\Bbb{Z}_2\times\Bbb{Z}_2\times\Bbb{Z}_2$ is the set of transformation $(x,y,z)\mapsto(\pm x,\pm y, \pm z)$. Composing $\tilde{f}$ with the degree four Belyi function $X(4)\rightarrow X(2)\cong\Bbb{CP}^1$ derived in  \eqref{Gamma(4)Belyi} shows that the desired  Belyi function  $f:X(8)\rightarrow X(2)\cong\Bbb{CP}^1$ of degree 32 is isomorphic to:
$$ X(8)\stackrel{\text{normalization}}{\dashrightarrow} \left\{(x,y,z)\in\Bbb{C}^3\mid x^2=z^8-1,y^2=z^4-1\right\}\stackrel{(x,y,z)\mapsto -\frac{1}{4}\frac{(z^4-1)^2}{z^4}}{\longrightarrow}\Bbb{CP}^1$$
The  group ${\rm{Deck}}\left(f:X\left(8\right)\rightarrow X\left(2\right)\right)\cong\bar{\Gamma}(2)/\bar{\Gamma}(8)$ is a central extension of $\bar{\Gamma}(4)/\bar{\Gamma}(8)\cong \Bbb{Z}_2\times\Bbb{Z}_2\times\Bbb{Z}_2$ by the group $\bar{\Gamma}(2)/\bar{\Gamma}(4)\cong \Bbb{Z}_2\times\Bbb{Z}_2$. According to 4.40 in \cite{dessin},  the group of automorphisms of this dessin   is the centralizer in ${\rm{S}}_{32}$ of the subgroup of permutations of elements of $\bar{\Gamma}(2)/\bar{\Gamma}(8)$ generated by right multiplication with the cosets of $\begin{bmatrix}
1&2\\
0&1
\end{bmatrix}, \begin{bmatrix}
1&0\\
-2&1
\end{bmatrix}$.
\end{example}

\begin{example}\label{noncongruence}
The simple description in \ref{secondBelyi} of the monodromy representation allows one to demonstrate examples of Belyi functions associated with non-congruence subgroups by exhibiting their permutation representations. For instance, define a surjective homomorphism from the free group on two generators $\bar{\Gamma}(2)$ onto ${\rm{A}_6}$ via sending  $T^2$ , $ST^2S^{-1}$ to the generators $(123)$, $(23456)$ of the alternating group. The kernel, denoted with $\Gamma$, is a normal subgroup of $\Gamma(2)$ of index $360$ which is not congruence since ${\rm{A}_6}$ is not a quotient of $\Gamma(2)/\Gamma(2N)$ for any $N$, cf. \cite{konrad}. Note that under $\Gamma(2)\rightarrow {\rm{A}_6}$:
$$T^2\mapsto (123)\quad  ST^2S^{-1}\mapsto (23456)\quad (TS)T^2(TS)^{-1}\mapsto (23456)^{-1}(123)^{-1}=(12)(6543)$$ 
Thus, following notations of \cite{dessin}, in the permutation representation $(\sigma_0,\sigma_1,\sigma_\infty)$ of dessin, $\sigma_0,\sigma_1,\sigma_\infty$ are the permutations of elements of ${\rm{A}}_6$ induced by right multiplication with $(23456)$, $(12)(6543)$, $(123)$, respectively. Plugging their orders $k_0=5$, $k_1=4$, $k_\infty=3$ and $m=360$ in the Corollary \ref{secondBelyisymmetric} yields $40$ for the genus of $X(\Gamma)$ and implies that the dessin  has $360$ edges, $72$ black vertices of degree $5$, $90$ white vertices of degree $4$ and $120$ faces each with $6$ edges. One can even extract some arithmetic information from this monodromy data: the group of permutations on elements of $\bar{\Gamma}(2)/\bar{\Gamma}\cong{\rm{A}}_6$ generated by the right action of a generating set of this quotient group obviously has the center equal to that of ${\rm{A}}_6$ which is trivial. Now the Corollary 3.2 from \cite{Debes-Douai} implies that the field of moduli for this Belyi function is a field of definition. More generally and by the same argument, this is the case for any Belyi function $X(\Gamma)\rightarrow X(2)$ where the group ${\rm{N}}_{\bar{\Gamma}(2)}(\bar{\Gamma})/\bar{\Gamma}$ is centerless. 
\end{example}

\section{Dessins on $X_0(N)$ and Modular Equation}
The primary objective of this section is to study the Belyi function $f:X_0(N)=\Gamma_0(N)\backslash\Bbb{H}^*\rightarrow X(1)=\Gamma(1)\backslash\Bbb{H}^*\cong\Bbb{CP}^1$, which according to Proposition \ref{firstBelyi} is actually $\frac{1}{1728}{j}$,  and employ it to compute the  {\it{ modular equation}}  associated with $\Gamma_0(N)$ which is the algebraic dependence relation between $\Gamma_0(N)$-invariant functions $z\mapsto j(z)$ and  $z\mapsto j(Nz)$. It is well-known that these two modular functions generate the function field of $X_0(N)$ and the minimal polynomial of $j(Nz)$ over $\Bbb{C}\left(j(z)\right)$ is with integer coefficients and of degree $\mu$ given by 
$$\mu:=\left[\bar{\Gamma}(1):\bar{\Gamma}_0(N)\right]=\left[\Gamma(1):\Gamma_0(N)\right]=N\prod_{\substack{p\mid N\\ p \text{ prime} }}(1+\frac{1}{p})$$  
See \cite{Milne} for further details. The meromorphic function $[z]\mapsto j(z)$ on the compact Riemann surface $X_0(N)$ may be determined provided that Belyi function $f:X_0(N)\rightarrow X(1)\cong\Bbb{CP}^1$ has been computed. Since the subgroup $\Gamma_0(N)$ is not normal,  this Belyi function is not regular and thus computing it is much harder than what carried out in section 2 examples due to lack of symmetry. But as $\langle T\rangle\Gamma(N)\leq\Gamma_0(N)$, once the regular Belyi function $X(N)\rightarrow X(1)$ and its very symmetric dessin are  obtained along with an explicit description of the action of the group ${\rm{Deck}}\left(X(N)\rightarrow X(1)\right)\cong\bar{\Gamma}(1)/\bar{\Gamma}(N)$ generated by automorphisms of $X(N)=\Gamma(N)\backslash\Bbb{H}^*$ obtained from left actions of 
$T=\begin{bmatrix}
1 & 1\\
0 & 1
\end{bmatrix}$ and
$S=\begin{bmatrix}
0 & -1\\
1 & 0
\end{bmatrix}$, one may factor this Belyi function and its dessin through the action of the subgroup $\bar{\Gamma}_0(N)/\bar{\Gamma}(N)$ of deck transformations in order to derive the Belyi function $X_0(N)\rightarrow X(1)$ and its dessin. We will illustrate this idea for  $N\in\{2,3,6\}$ where $\bar{\Gamma}_0(N)$ and $\langle T\rangle\bar{\Gamma}(N)$ coincide. So to obtain $X_0(N)\rightarrow X(1)$ it only suffices to consider the quotient of $X(N)\rightarrow X(1)$  under the order $N$ automorphism induced by $T$.

\begin{example}\label{Gamma_0(2)}
Suppose $N=2$. In Example \ref{Gamma(2)} the Belyi function $X(2)\cong\Bbb{CP}^1\rightarrow  X(1)\cong\Bbb{CP}^1$  was calculated as $f(z)=\frac{(3z^2+1)^3}{(9z^2-1)^2}$ (cf. \eqref{Gamma(2)Belyi}) whose dessin is depicted in Figure 1 and its group of deck transformations is a version of ${\rm{S}}_3$ generated by involutions $z\mapsto -z$ and $z\mapsto\frac{-z+1}{3z+1}$ corresponding to $S$ and $T$, respectively. The latter is an order two M\"obius transformation with fixed points $-1,\frac{1}{3}$. Therefore, any meromorphic function $\Bbb{CP}^1\rightarrow\Bbb{CP}^1$ invariant under it factors via $z\mapsto\left(\frac{z+1}{z-\frac{1}{3}}\right)^2$. Applying this map to coordinates of vertices in Figure 1, we arrive at Figure 6.
\begin{figure}[ht]
\centering
\includegraphics[height=4.5cm]{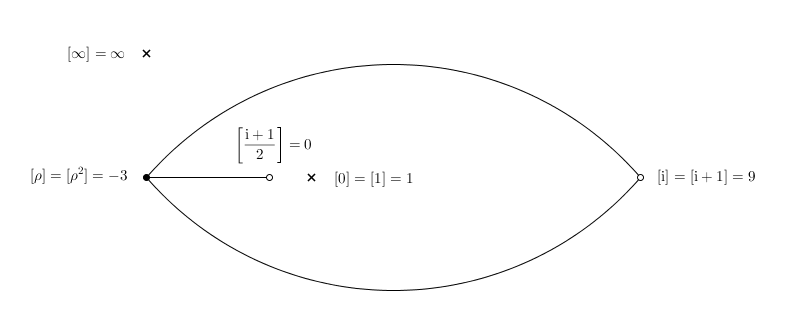}\\
Figure 6
\end{figure}
 \noindent Note that some of the orbits of the action of $\Gamma(2)$ appeared in Figure 1 are identified in Figure 6 under the action of $T$  on the upper half plane given by $z\mapsto z+1$. Belyi function of this dessin is $\frac{1}{27}\frac{(z+3)^3}{(z-1)^2}$ and satisfies $\frac{1}{27}\frac{(z+3)^3}{(z-1)^2}-1=\frac{1}{27}\frac{z(z-9)^2}{(z-1)^2}$. In summary:
\begin{equation*}%\label{Gamma_0(2)Belyi}
\xymatrixcolsep{8pc}\xymatrix{X(2)=\Gamma(2)\backslash\Bbb{H}^*\cong\Bbb{CP}^1 \ar[d]_{z\mapsto\left(\frac{z+1}{z-\frac{1}{3}}\right)^2 }  \ar[dr]^{z\mapsto\frac{(3z^2+1)^3}{(9z^2-1)^2}}\\ 
X_0(2)=\Gamma_0(2)\backslash\Bbb{H}^*\cong\Bbb{CP}^1 \ar[r]_{z\mapsto\frac{1}{27}\frac{(z+3)^3}{(z-1)^2}} & X(1)=\Gamma(1)\backslash\Bbb{H}^*\cong\Bbb{CP}^1 }
\end{equation*}
 According to Proposition \ref{firstBelyi}, the dessin in Figure 6 conveys information about the congruence subgroup $\Gamma_0(2)$: it possesses exactly one elliptic orbit, that of $\frac{{\rm{i}}+1}{2}$\footnote{The matrix $\begin{bmatrix} 
-1& 1\\
-2&1
\end{bmatrix}$ in $\Gamma_0(2)$ fixes $\frac{{\rm{i}}+1}{2}$.}
 which is of order two \footnote{ In general, $\Gamma_0(N)$ has an elliptic point of order three iff $N$ is odd and $-3$ is a quadratic residue modulo $N$ and has an elliptic point of order two iff $-1$  is a quadratic residue modulo $N$.} . Looking at  $\times$ vertices, there are two cusps $[0],[\infty]$ of degrees two and one respectively\footnote{The cusp $\infty$ of $\Gamma_0(N)$ is always of width  one as its stabilizer $\langle T\rangle$ is contained in $\Gamma_0(N)$ while the width of cusp $0$ is $N$.}.  
\end{example}
\begin{example}\label{Gamma_0(3)}
Suppose $N=3$. $\bar{\Gamma}_0(3)=\langle T\rangle\bar{\Gamma}(3)$ is an index $N+1=4$ subgroup of $\bar{\Gamma}(1)={\rm{PSL}}_2(\Bbb{Z})$. So the degree twelve Belyi function 
$\begin{cases}
X(3)\rightarrow X(1)\\
z\mapsto \frac{1}{64}\frac{z^3(z^3+8)^3}{(z^3-1)^3}
\end{cases}$
 written in \eqref{Gamma(3)Belyi} factors  in the tower below of meromorphic functions on $\Bbb{CP}^1$: 
$$X(3)=\Gamma(3)\backslash\Bbb{H}^*\stackrel{3-\text{sheeted}}{\longrightarrow}X_0(3)=\Gamma_0(3)\backslash\Bbb{H}^*\stackrel{4-\text{sheeted}}{\longrightarrow}X(1)=\Gamma(1)\backslash\Bbb{H}^*$$
The first map is the quotient map for the action of the group $\langle T\rangle$ on $\Bbb{CP}^1$. We saw in \ref{Gamma(3)} that $T$ acts as the rotation  through $120\,^\circ$ about the origin.  Hence, this quotient map can be considered as $z\mapsto z^3$ and we arrive at $z\mapsto\frac{1}{64}\frac{z(z+8)^3}{(z-1)^3}$ as the equation of the Belyi function $X_0(3)\cong\Bbb{CP}^1\rightarrow X(1)\cong\Bbb{CP}^1$ whose dessin is illustrated in Figure 7 and is  obtained from applying $z\in\Bbb{CP}^1\mapsto z^3\in\Bbb{CP}^1$ (induced by $T:z\in\Bbb{H}\rightarrow z+1\in\Bbb{H}$) to the vertices in Figure 3.
\begin{figure}[ht]
\centering
\includegraphics[height=4.5cm]{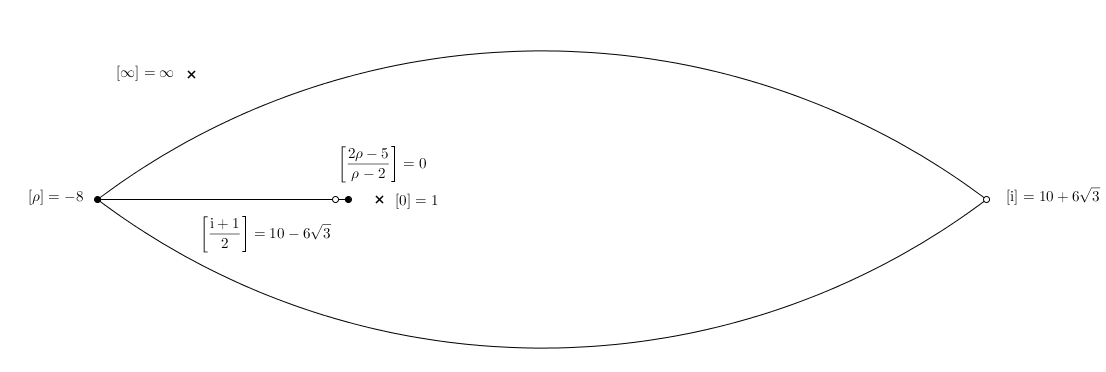}\\
Figure 7
\end{figure}
\noindent Note that the Belyi property translates to $\frac{1}{64}\frac{z(z+8)^3}{(z-1)^3}-1=\frac{1}{64}\frac{(z^2-20z-8)^2}{(z-1)^3}$. Again,  some data concerning the group  $\Gamma_0(3)$ might be read off from the dessin in Figure 7:  $\left[\frac{2\rho-5}{\rho-2}\right]$ is the unique elliptic orbit whose order is three whereas having a degree one black vertex, it precisely has one elliptic orbit of order three\footnote{The matrix 
$\begin{bmatrix}
-8&19\\
-3&7
\end{bmatrix}$
in $\Gamma_0(3)$ fixes $\frac{2\rho-5}{\rho-2}$ where as usual $\rho={\rm{e}}^{\frac{\pi{\rm{i}}}{3}}$.}. Analyzing  $\times$ vertices implies that orbits of cusps for  $\Gamma_0(3)$ are $[\infty]$ and $[0]$ of widths one and three respectively. 
\end{example}

\begin{example}\label{Gamma_0(6)}
For $N=6$, we have equalities $\mu_6=\left[\bar{\Gamma}(1):\bar{\Gamma}(6)\right]=72,\left[\bar{\Gamma}(1):\bar{\Gamma}_0(6)\right]=6\left(1+\frac{1}{2}\right)\left(1+\frac{1}{3}\right)=12$ while  the coset of $T$  in $\bar{\Gamma}_0(6)/\bar{\Gamma}(6)$ is of order $6=\frac{72}{12}$. We deduce that $\bar{\Gamma}_0(6)$ coincides with $\langle T\rangle\bar{\Gamma}(6)$ and therefore the degree twelve Belyi function $X_0(6)\rightarrow X(1)\cong\Bbb{CP}^1$ is just  the quotient of  map $X(6)\rightarrow X(1)\cong\Bbb{CP}^1$ appeared in \eqref{Gamma(6)Belyi} under the automorphism induced by the action of $T$  which is an order six deck transformation. But $X(6)$ is the hexagonal elliptic curve $y^2=x^3-1$ whose group of automorphisms is of order six generated by $(x,y)\mapsto\left({\rm{e}}^{\frac{2\pi{\rm{i}}}{3}}x,-y\right)$. This indicates that one just needs to form the quotient of this elliptic curve under this automorphism which is:
$\begin{cases}
\{y^2=x^3-1\}\rightarrow\Bbb{CP}^1\\
(x,y)\mapsto x^3
\end{cases}$.
Therefore, the final answer is the degree twelve Belyi function on $\Bbb{CP}^1$ obtained by replacing $x^3$ with $x$ in right-hand side of \eqref{Gamma(6)Belyi} (as usual, coordinate on the Riemann sphere is denoted by $z$):
\begin{equation}\label{Gamma_0(6)Belyi}
z\mapsto -\frac{(z-4)^3(z^3-228z^2+48z-64)^3}{1728z^2(z-1)^3(z+8)^6}
\end{equation}
where (using the second identity in \eqref{Gamma(3)Belyi} to obtain its difference with $1$) the Belyi property is reflected in the identity:
$$ -\frac{(z-4)^3(z^3-228z^2+48z-64)^3}{1728z^2(z-1)^3(z+8)^6}-1=-\frac{\left(\left(z^3+258z^2+48z-64\right)^2-78732z^4\right)^2}{1728z^2(z-1)^3(z+8)^6}$$
Note that as $2,3\mid N=6$, $\Gamma_0(N)$ does not have any elliptic element and thus from Proposition \ref{firstBelyi}, its twelve edge dessin has four black vertices of degree three and six white vertices of degree two. We infer that the underlying graph is the same as that of Example \ref{Gamma(3)}. Nevertheless, the Belyi functions written in \eqref{Gamma(3)Belyi} and \eqref{Gamma_0(6)Belyi} are not isomorphic since the former is regular and the latter is not, cf. Remark \ref{complexstructure}.
\end{example}

\begin{example}\label{Gamma_0(11)}
Here we depict the dessin corresponding to the congruence subgroup $\Gamma_0(11)$. This subgroup is of index  12, does not possess any elliptic element and has exactly two cusps, namely the orbits of $0,\infty$, of widths 11,1. Hence, some basic properties of its dessin can be inferred: it is a dessin on the torus with 12 edges, 4 black vertices of degree three, 6 white vertices of degree two and 2 faces surrounded by 22 and 2 edges. Therefore, by deleting white vertices, we have an embedding of a graph with four degree three vertices that has exactly one loop (which is the boundary of the face with two edges and contributes two to the degree of the corresponding vertex) on the torus  as a map whose complement consists of two $2$-cells where of course one of them is the one determined by the loop. The reader may quickly check that essentially there is only one way to draw such a graph (as a map) on the torus and by putting back one white vertex on each edge, we arrive at the dessin of $X_0(11)$  in Figure 8.
\begin{figure}[ht]
\centering
\includegraphics[height=5cm]{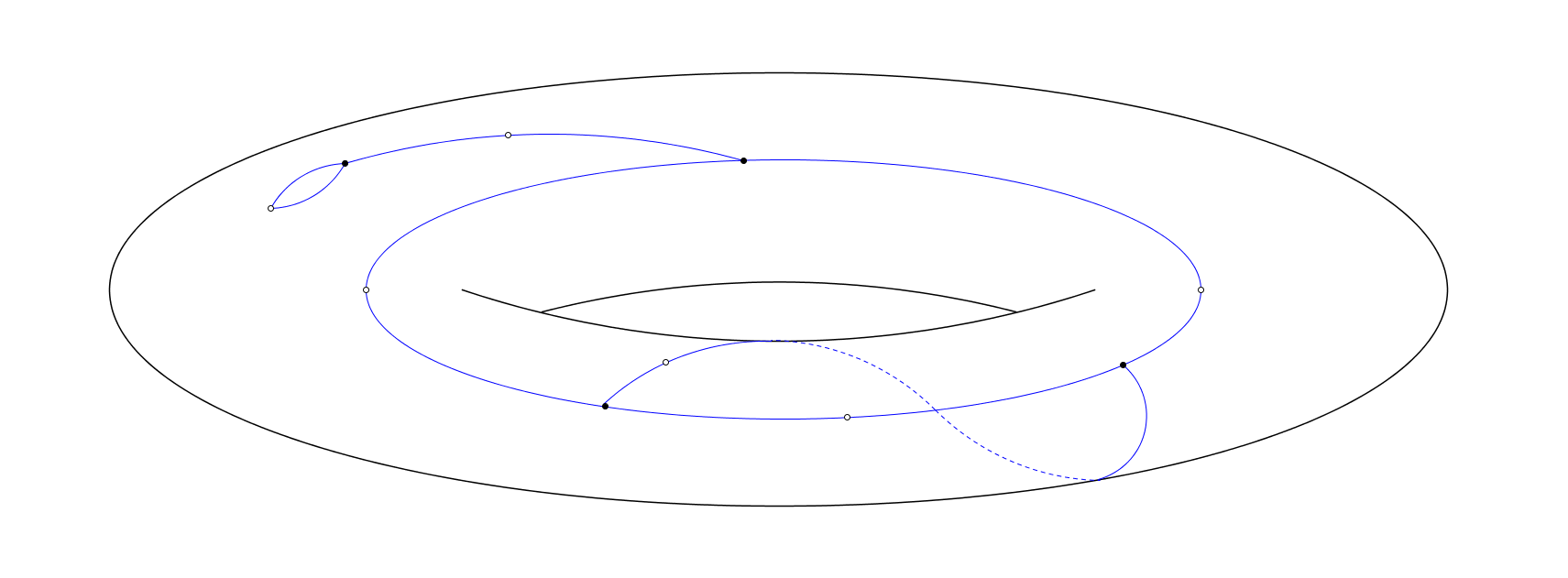}\\
Figure 8
\end{figure}
\end{example}

\noindent After studying the Belyi function $f:[z]\in X_0(N)\mapsto\frac{1}{1728}j(z)$ in these examples, let us concentrate on $z\in\Bbb{H}\mapsto j(Nz)$ which is $\Gamma_0(N)$-invariant and therefore gives a meromorphic function on $X_0(N)$. Once more, this function has at most three ramification values, i.e. $j(\infty)=\infty$ and  the ramification values $j(\rho)=0$ and $j({\rm{i}})=1728$ of $j:\Bbb{H}\rightarrow\Bbb{C}$. Thus, normalizing it, one gets a new Belyi function on $X_0(N)$ denoted by $g$.  Calculating first equations for $f$ and $g$ and then the algebraic dependence relation between them  results in the modular equation for $\Gamma_0(N)$. This is our key idea: thinking of the modular equation for $\Gamma_0(N)$ as an algebraic relation that two Belyi functions $f:[z]\mapsto\frac{1}{1728}j(z)$ and $g:[z]\mapsto\frac{1}{1728}j(Nz)$ on the modular curve $X_0(N)$ should satisfy.\\
\begin{remark}\label{rigid}
It is vital to note that in our approach to computing the modular equation it is not sufficient to only obtain Belyi functions $f$ and $g$ up to isomorphism since we are interested in the algebraic dependence relation between them. Thus, after computing a version of the Belyi function $f:X_0(N)\rightarrow\Bbb{CP}^1$ (corresponding to $j(z)$), the task that has been carried out in Examples \ref{Gamma_0(2)} , \ref{Gamma_0(3)} and \ref{Gamma_0(6)} for $N\in\{2,3,6\}$, its dessin is fixed in  both algebraic (in terms of the $\Gamma_0(N)$-orbits  appearing as vertices) and geometric (in the sense of coordinates on the Riemann surface $X_0(N)$) sense and then one should consider how the dessin of the Belyi function $g$ (corresponding to $j(Nz)$) is positioned  with respect to it. For $N=2,3$ this is going to be achieved in Examples \ref{modular2} , \ref{modular3}.    
\end{remark}

 What can be said about $g$? Black (resp. white) vertices of its dessin correspond to $\Gamma_0(N)$-orbits of points $z$ in the upper half plane lying in the $\Gamma_0(N)$-invariant subset $\frac{1}{N}\Gamma(1).\rho$ (resp. $\frac{1}{N}\Gamma(1).{\rm{i}}$) where again just as the dessin of $f$ described in Proposition \ref{firstBelyi}, their degree is either one or three (resp. two) where degree one occurs precisely when $z$ is an elliptic point of $\Gamma_0(N)$. Just like the case of $f$, vertices of type $\times$ (poles of the Belyi function $g$) are orbits of cusps for $\Gamma_0(N)$. The degree of  $g$ is same as that of $f$, i.e. $\mu=\left[\Gamma(1):\Gamma_0(N)\right]=\left[\bar{\Gamma}(1):\bar{\Gamma}_0(N)\right]$. These are immediate consequences of the following commutative diagram: 
$$\xymatrix{\Bbb{H} \ar[d]   \ar[dr]^{[z]\mapsto j(Nz)}& \\
                 \Gamma_0(N)\backslash\Bbb{H}\ar[r] & \Gamma(1)\backslash\Bbb{H}=X(1)-[\infty]\cong\Bbb{C}}$$

\noindent With some group-theoretic manipulations, one can determine the orbits that are vertices of the dessin of $g$  on $X_0(N)$: choose a complete set of left cosets  representatives of $\Gamma_0(N)$ in $\Gamma(1)$ denoted by $\gamma_1,\dots,\gamma_\mu$. It is easy to verify that:
\begin{equation}\label{Gamma_0(N)vertex}
 \forall\tau\in\Bbb{H}: \frac{1}{N}\Gamma(1).\tau=\bigcup_{j=1}^\mu\Gamma_0(N).\left(\frac{\gamma_j^{\rm{t}}.\tau}{N}\right)
\end{equation}
As the final remark that will aid us below in drawing  the dessins of $f$ and $g$, consider points $[0],[\infty]\in X_0(N)=\Gamma_0(N)\backslash\Bbb{H}^*$  which are vertices of type $\times$ in both cases. According to \ref{firstBelyi}, their degrees in the dessin of $f$    are given by $N$ and $1$ respectively. These degrees are reversed in the dessin of $g$: writing down the $q$-expansion and carefully analyzing charts around these points of $X_0(N)$, it is not hard to see that the orders of poles at $[0]$ and $[\infty]$ of $[z]\in X_0(N)\mapsto j(Nz)$ are $1$ and $N$ respectively.  \\
\begin{example}\label{modular2}
Let us study the Belyi function $g:[z]\mapsto \frac{1}{1728}j(2z)$ on $X_0(2)$, where  
 a complete set of left  coset representatives of $\Gamma_0(2)$ in $\Gamma(1)$ is given by:
$$\left\{{\rm{I}}_2,TS=\begin{bmatrix}
1&-1\\
1&0
\end{bmatrix}, (TS)^2=\begin{bmatrix}
0&-1\\
1&-1
\end{bmatrix}\right\}$$
 Action of their transpose matrices is given by $(TS)^{\rm{t}}:z\mapsto -1-\frac{1}{z}$ and $\left((TS)^2\right)^{\rm{t}}:z\mapsto \frac{-1}{z+1}$. Now, using \eqref{Gamma_0(N)vertex}:
\begin{equation} \label{Gamma_0(2)vertex}
\begin{split}
&\frac{1}{2}\Gamma(1).\rho=\Gamma_0(2).\left(\frac{\rho}{2}\right) \bigcup \Gamma_0(2).\left(\frac{\rho}{2}-1\right)\bigcup \Gamma_0(2).\left(\frac{\rho-2}{6}\right) \\
& \frac{1}{2}\Gamma(1).{\rm{i}}=\Gamma_0(2).\left(\frac{{\rm{i}}}{2}\right) \bigcup \Gamma_0(2).\left(\frac{{\rm{i}}-1}{2}\right)\bigcup \Gamma_0(2).\left(\frac{{\rm{i}}-1}{4}\right)
\end{split}
\end{equation}
Orbits of $\Gamma_0(2)$ appeared above are black and white vertices of the dessin  on $X_0(2)$ corresponding to $g$. But we should omit repetitions of orbits: 
$\frac{\rho}{2},\frac{\rho}{2}-1,\frac{\rho-2}{6}$ are congruent modulo the  action of $\Gamma_0(2)$: $\begin{bmatrix}
1&-1\\
0&1
\end{bmatrix},\begin{bmatrix}
1&0\\
-4&1
\end{bmatrix}\in\Gamma_0(2)$ send the first point to the other two. Hence, there is only one black vertex, $\left[\frac{\rho}{2}\right]$ whose degree is three as $\Gamma_0(2)$ does not have any elliptic element of the order three. For  white vertices, in the second line of \eqref{Gamma_0(2)vertex} $\frac{{\rm{i}}}{2}$ and $\frac{{\rm{i}}-1}{2}$ are inequivalent under the action of $\Gamma(1)$ while $\begin{bmatrix}
1&0\\
-2&1
\end{bmatrix}\in\Gamma_0(2)$ maps $\frac{{\rm{i}}}{2}$ to $\frac{{\rm{i}}-1}{4}$. We deduce that there are two white vertices: $\left[\frac{{\rm{i}}}{2}\right]$ and $\left[\frac{{\rm{i}}-1}{2}\right]=\left[\frac{{\rm{i}}+1}{2}\right]$. The degree of such a vertex is one iff it is the unique order two elliptic orbit of $\Gamma_0(2)$ and otherwise equals two. According to the dessin in  Figure 6 and the discussion after it, $\left[\frac{{\rm{i}}+1}{2}\right]$ is the unique elliptic orbit of $\Gamma_0(2)$. Up to now, the dessin of $g$  on $X_0(2)\cong\Bbb{CP}^1$  possesses three edges, one black vertex of degree three and two white vertices of degrees one and two which implies that there are two faces. But as we observed before, points $[0],[\infty]$ are two vertices of degrees $1$ and $N=2$ each corresponds to a face. The coordinates of $\left[\frac{{\rm{i}}+1}{2}\right],[0],[\infty]$ may be read off from the dessin of $f$  in Figure 6 which implies that $g(z)$ must be in the form of $\frac{k(z-\alpha)^3}{z-1}$  satisfying $g(z)-1=\frac{kz(z-\beta)^2}{z-1}$ where $\alpha=\left[\frac{\rho}{2}\right], \beta=\left[\frac{{\rm{i}}}{2}\right]$. Fixing three vertices of a dessin on $\Bbb{CP}^1$ rigidifies it so the unknowns $\alpha,\beta$ and $k$ are readily determined: $k=\frac{64}{27},\alpha=\frac{3}{4},\beta=\frac{9}{8}$. Thus
 $$g(z)=\frac{1}{27}\frac{(4z-3)^3}{z-1}\quad g(z)-1=\frac{1}{27}\frac{z(8z-9)^2}{z-1}$$ 
and the dessin is depicted in Figure 9.

\begin{figure}[ht]
\centering
\includegraphics[height=4.5cm]{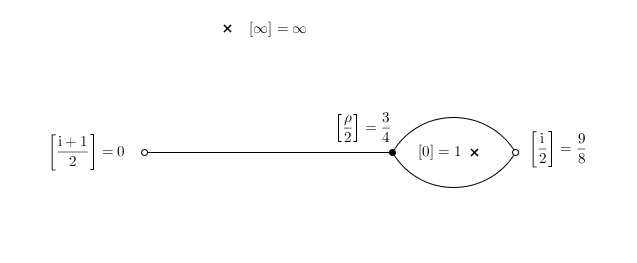}\\
Figure 9
\end{figure}

\noindent Combining this with the results of Example \ref{Gamma_0(2)}:
\begin{equation}\label{pair1}
f(z)=\frac{1}{27}\frac{(z+3)^3}{(z-1)^2}\quad  g(z)=\frac{1}{27}\frac{(4z-3)^3}{z-1}\
\end{equation}
whose dessins are illustrated in Figures 6 and 9. Therefore, we arrive at equations $1728f=64\frac{(z+3)^3}{(z-1)^2}$ and $1728g=64\frac{(4z-3)^3}{z-1}$ for meromorphic functions $[z]\mapsto j(z)$ and $[z]\mapsto j(2z)$ on $X_0(2)=\Gamma_0(2)\backslash\Bbb{H}^*\cong\Bbb{CP}^1$ and we must find their algebraic dependence relation. A simple observation helps us to accomplish this:  changing $z$ to $z+1$ in them yields  new rational functions $64\frac{(z+4)^3}{z^2}$ and     $64\frac{(4z+1)^3}{z}$ where the former transforms to the latter by substituting $z$ with $\frac{1}{z}$. We deduce that any polynomial equation that $1728f$ and $1728g$ satisfy is symmetric in them and actually a polynomial relation between $64\frac{(z+4)^3}{z^2}+64\frac{(4z+1)^3}{z}$ and $\left(64\frac{(z+4)^3}{z^2}\right)\left(64\frac{(4z+1)^3}{z}\right)$ which are polynomials in $y:=z+\frac{1}{z}$. 
\begin{equation*}
\begin{split}
&64\frac{(z+4)^3}{z^2}+64\frac{(4z+1)^3}{z}=4096\left(z^2+\frac{1}{z^2}\right)+3136\left(z+\frac{1}{z}\right)+1536= 4096y^2+3136y-6656 \\
&\left(64\frac{(z+4)^3}{z^2}\right)\left(64\frac{(4z+1)^3}{z}\right)=  262144\left(z^3+\frac{1}{z^3}\right)+3342336\left(z^2+\frac{1}{z^2}\right)+14991360\left(z+\frac{1}{z}\right)\\
&+26808320=262144y^3+3342336y^2+14204928y+20123648
\end{split}
\end{equation*}
One can compute the algebraic dependence relation between the above quadratic and cubic in $y$ (which is straightforward using the division algorithm) and then in the derived equation replace them with $X+Y$ and $XY$ respectively. The result is  the modular equation for $\Gamma_0(2)$:
\begin{equation*}
\begin{split}
 &X^3+Y^3-X^2Y^2+1488(X^2Y+XY^2)-162000
(X^2+Y^2)\\
&+40773375XY+8748000000(X+Y)-157464000000000=0
\end{split}
\end{equation*}
\end{example}

\begin{example}\label{modular3}
Here, we derive  an equation for the  Belyi function $g:[z]\mapsto \frac{1}{1728}j(3z)$ on $X_0(3)\cong\Bbb{CP}^1$.  Let us work with the following complete set of representatives:
$$\{\gamma_1,\gamma_2,\gamma_3,\gamma_4\}=\left\{{\rm{I}}_2,S,TST,STST=T^{-1}S\right\}$$
wherein $\bar{\Gamma}(1)={\rm{PSL}}_2(\Bbb{Z})<{\rm{Aut}}(\Bbb{H})$  their transpose matrices may be written as:
$${\rm{I}}_2:z\mapsto z\quad S^{\rm{t}}=S: z\mapsto\frac{-1}{z}\quad (TST)^{\rm{t}}=T:z\mapsto z+1\quad (T^{-1}S)^{\rm{t}}=TS:z\mapsto\frac{z-1}{z}$$ 
Substituting in \eqref{Gamma_0(N)vertex} and noting that $S$ and $TS$ stabilize ${\rm{i}}$ and $\rho$ respectively:
\begin{equation} \label{Gamma_0(3)vertex}
\begin{split}
&\frac{1}{3}\Gamma(1).\rho=\Gamma_0(3).\left(\frac{\rho}{3}\right) \bigcup \Gamma_0(3).\left(\frac{\rho-1}{3}\right)\bigcup \Gamma_0(3).\left(\frac{\rho+1}{3}\right)\\
& \frac{1}{3}\Gamma(1).{\rm{i}}=\Gamma_0(3).\left(\frac{{\rm{i}}}{3}\right) \bigcup \Gamma_0(3).\left(\frac{{\rm{i}}+1}{3}\right)
\end{split}
\end{equation}
We should exclude repeated   $\Gamma_0(3)$-orbits in \eqref{Gamma_0(3)vertex} in order to get black and white vertices of the dessin  on $X_0(3)$ corresponding to $g$. Orbits  $\left[\frac{\rho}{3}\right]$ , $\left[\frac{\rho+1}{3}\right]$  coincide with $\left[\frac{2\rho-5}{\rho-2}\right]$  appeared as a vertex of the dessin in Figure 7 ($T^2(\frac{\rho+1}{3})=\frac{2\rho-5}{\rho-2}$ and $ST^3S\left(\frac{\rho}{3}\right)=\frac{\rho-1}{3}$ where $T^2,ST^3S\in\Gamma_0(3)$) while it can be easily verified that no element of $\Gamma(1)$ carries $\frac{\rho}{3}$ to $\frac{\rho+1}{3}$. Hence, there are two black vertices: $\left[\frac{\rho}{3}\right],\left[\frac{2\rho-5}{\rho-2}\right]$. The degree of a black vertex is one iff it is an order three elliptic orbit of $\Gamma_0(3)$. Following discussions in Example \ref{Gamma_0(3)} based on the  dessin in Figure 7, $\left[\frac{2\rho-5}{\rho-2}\right]$  is the only elliptic orbit of $\Gamma_0(3)$.  For  white vertices, in the second line of \eqref{Gamma_0(3)vertex} $\frac{{\rm{i}}}{3}$ and $\frac{{\rm{i}}+1}{2}$ are inequivalent under the action of $\Gamma(1)$. We deduce that there are two white vertices: $\left[\frac{{\rm{i}}}{3}\right]$ and $\left[\frac{{\rm{i}}+1}{3}\right]$,  each one of degree  two  because $\Gamma_0(3)$ does not have any elliptic point of order two. There are two faces  $[0],[\infty]$ which form the set of cusps' orbits of $\Gamma_0(3)$ and  their degrees are $1$ and $N=3$ respectively. This gives us the depiction of the dessin of $g$  in Figure 10  where the coordinates of $\left[\frac{2\rho-5}{\rho-2}\right],[0],[\infty]$ have been obtained from the dessin of $f$   (Figure 7)\footnote{The same argument as that of the Remark \ref{draw} shows that  in the dessin of the  Belyi function $[z]\in X_0(N)\mapsto\frac{1}{1728}j(Nz)$ the edges are always in the form of  $\left\{\left[\frac{1}{N}\gamma.{\rm{e}}^{{\rm{i}}\theta}\right] \big | \frac{\pi}{3}\leq \theta\leq\frac{\pi}{2}\right\}$ where $\gamma\in\Gamma_0(N)$. When $N=3$, this fact guides us in drawing  the dessin in Figure 10: vertices corresponding to orbits of $\frac{1}{3}T(\rho)$ and $\frac{1}{3}T({\rm{i}})=\left[\frac{{\rm{i}}+1}{3}\right]$ should be connected. The former is $\left[\frac{2\rho-5}{\rho-2}\right]$ as $T^2\in\Gamma_0(3)$ maps  $\frac{1}{3}T(\rho)=\frac{\rho+1}{3}$ to $\frac{2\rho-5}{\rho-2}$.}. Again, these three coordinates rigidify the dessin and the coordinates of unknown vertices $\left[\frac{{\rm{i}}}{3}\right]$, $\left[\frac{{\rm{i+1}}}{3}\right]$ and $\left[\frac{\rho}{3}\right]$ will be uniquely determined. We leave the details to the reader and just exhibit the final picture in Figure 10 and an equation for $g$  after that:

\begin{figure}[ht]
\centering
\includegraphics[height=4.5cm]{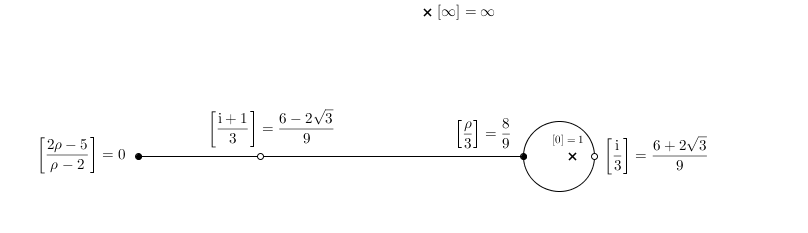}\\
Figure 10
\end{figure}
$$g(z)=\frac{1}{64}\frac{z(9z-8)^3}{z-1} \quad g(z)-1=\frac{1}{64}\frac{(27z^2-36z+8)^2}{z-1} $$
Recalling results of Example \ref{Gamma_0(3)}, our pair of Belyi functions on $X_0(3)\cong\Bbb{CP}^1$ is:
\begin{equation}\label{pair2}
f(z)=\frac{1}{64}\frac{z(z+8)^3}{(z-1)^3}\quad g(z)=\frac{1}{64}\frac{z(9z-8)^3}{z-1}
\end{equation}
Compare their dessins in Figures 7 and 10. Next task is to compute the algebraic equation that meromorphic functions $1728f=27\frac{z(z+8)^3}{(z-1)^3}$ and $1728g=27\frac{z(9z-8)^3}{z-1}$ satisfy. Again, the pair $(f,g)$ has a property that facilitates this: $f\left(\frac{z}{z-1}\right)=g(z)$, i.e. after substituting $z$ with $z+1$ in  $27\frac{z(z+8)^3}{(z-1)^3}$ and $27\frac{z(9z-8)^3}{z-1}$, the resulting expressions $27\frac{(z+1)(z+9)^3}{z^3}$ and $27\frac{(z+1)(9z+1)^3}{z}$ transform to each other under $z\rightsquigarrow\frac{1}{z}$. Consequently, everything reduces to finding the algebraic dependence  relation between sum and product of two rational functions that were just described which themselves are polynomials in $y:=z+\frac{1}{z}$:
\begin{equation*}
\begin{split}
&27\frac{(z+1)(z+9)^3}{z^3}+27\frac{(z+1)(9z+1)^3}{z}=19683\left(z^3+\frac{1}{z^3}\right)+26244\left(z^2+\frac{1}{z^2}\right)\\
&+7317\left(z+\frac{1}{z}\right)+1512=19683y^3+26244y^2-51732y-50976\\
&\left(27\frac{(z+1)(z+9)^3}{z^3}\right)\left(27\frac{(z+1)(9z+1)^3}{z}\right)= 531441\left(z^4+\frac{1}{z^4}\right)+15588936\left(z^3+\frac{1}{z^3}\right)\\
&+163526364\left(z^2+\frac{1}{z^2}\right)+713411064\left(z+\frac{1}{z}\right)+1129884390\\
&=531441y^4+15588936y^3+161400600y^2+666644256y+803894544
\end{split}
\end{equation*}
It only remains to derive the algebraic dependence relation between the cubic and the quartic in $y$ appeared above (this is a bit harder and  best  carried out with a computer algebra package such as MAPLE) and then substitute the cubic (which came from $f+g$) with $X+Y$ and the quartic (which came from $fg$) with $XY$. This culminates in the modular equation for  $\Gamma_0(3)$ below:
\begin{equation*}
\begin{split}
 &X^4+Y^4-X^3Y^3+2232(X^3Y^2+X^2Y^3)+2587918086X^2Y^2-1069956(X^3Y+XY^3)\\
&+36864000(X^3+Y^3)+8900222976000(X^2Y+XY^2)+452984832000000(X^2+Y^2)\\
&-770845966336000000XY+1855425871872000000000(X+Y)=0
\end{split}
\end{equation*}
\end{example}

The existence of a relation such as $f\left(\frac{z}{z-1}\right)=g(z)$ between two Belyi functions in either of \eqref{pair1} or \eqref{pair2} that  facilitated computing the algebraic dependence for either of these pairs is not accidental. $z\mapsto \frac{z}{z-1}$ is an involution of $X_i(2)\cong\Bbb{CP}^1\,(i\in\{2,3\})$  composition by which in either of pairs swaps the functions in that pair. An involution of $X_0(N)$ with this property always exists and is the well-defined holomorphic involution $[z]\mapsto \left[\frac{-1}{Nz}\right]$ of $X_0(N)=\Gamma_0(N)\backslash\Bbb{H}^*$.  Composing it with the Belyi function $f:[z]\mapsto\frac{1}{1728}j(z)$ yields $[z]\mapsto\frac{1}{1728}j\left(\frac{-1}{Nz}\right)=j(Nz)$ which is the Belyi function $g$. Since under this involution $f$ goes to $g$, the symmetry of all algebraic equations that $f,g$ satisfy and in particular the modular equation will readily be inferred.  We summarize our method:
\begin{prop}\label{modularequation}
For any $N\in\Bbb{N}$, there is a pair of Belyi functions consists of 
$$\begin{cases}
f:X_0(N)\rightarrow\Bbb{CP}^1\\
[z]\mapsto\frac{1}{1728}j(z)
\end{cases}\quad \begin{cases}
g:X_0(N)\rightarrow\Bbb{CP}^1\\
[z]\mapsto\frac{1}{1728}j(Nz)
\end{cases}$$  
on the compact Riemann surface $X_0(N)=\Gamma_0(N)\backslash\Bbb{H}^*$. They generate the function field of $X_0(N)$  and transform to each other after being composed with the holomorphic involution $[z]\mapsto \left[\frac{-1}{Nz}\right]$. The minimal polynomial of $1728g$ over $\Bbb{C}(f)$ is symmetric and precisely the modular equation for $\Gamma_0(N)$.
\end{prop}

We finish with another interesting consequences of the constructions used  in the last two examples. Suppose that in Proposition \ref{modularequation} after fixing the dessin of $f$  (cf. Remark \ref{rigid}), $\Gamma_0(N)$-orbits that occur as vertices of the dessin of $g$  have been determined after necessary group-theoretic arguments and also the positions of these vertices on $X_0(N)$ are known.  Evaluating $j=1728f$ at them might lead to some non-trivial $j$-invariant calculations because among the vertices of the dessin of $g$  there are some new orbits that had not appeared previously. For example, compare Figures 6 and 9 where $N=2$ and Figures 7 and 10 where $N=3$. In the first case, evaluating the Belyi function  $f(z)=\frac{1}{27}\frac{(z+3)^3}{(z-1)^2}$ derived in \eqref{pair1} at $\left[\frac{\rho}{2}\right]=\frac{3}{4}$ and   $\left[\frac{{\rm{i}}}{2}\right]=\frac{9}{8}$ leads to: 
$$j\left(\frac{\rho}{2}\right)=1728.f\left(\frac{3}{4}\right)=16\times 15^3 \quad j\left(\frac{{\rm{i}}}{2}\right)=1728.f\left(\frac{9}{8}\right)=66^3$$
and in the second case, in Example \ref{modular3} one should compute the values of  $f(z)=\frac{1}{64}\frac{z(z+8)^3}{(z-1)^3}$  at points $\left[\frac{\rho}{3}\right]=\frac{8}{9}$, $\left[\frac{{\rm{i}}+1}{3}\right]=\frac{6-2\sqrt{3}}{9}$ and $\left[\frac{{\rm{i}}}{3}\right]=\frac{6+2\sqrt{3}}{9}$:
\begin{equation*} 
\begin{split}
& j\left(\frac{\rho}{3}\right)=1728.f\left(\frac{8}{9}\right)=-3\times 160^3\\
& j\left(\frac{{\rm{i}}+1}{3}\right)=1728.f\left(\frac{6-2\sqrt{3}}{9}\right)=\left(18-6\sqrt{3}\right)(82-54\sqrt{3})^3\\
& j\left(\frac{{\rm{i}}}{3}\right)=1728.f\left(\frac{6+2\sqrt{3}}{9}\right)=\left(18+6\sqrt{3}\right)(82+54\sqrt{3})^3
\end{split}
\end{equation*}
Another interesting set of values of the modular $j$-function can be obtained by evaluating $1728f$ at fixed points of the involution mentioned in Proposition \ref{modularequation}. In both of Examples \ref{modular2} and \ref{modular3} this involution was given by $z\mapsto\frac{z}{z-1}$ whose fixed points are $z=0,2$. When $N=2$, under the action of $\Gamma_0(2)$ we have $[z]=[\frac{-1}{2z}]$ for  $z=\frac{\sqrt{2}{{\rm{i}}}}{2},z=\frac{{\rm{i}}-1}{2}$. Since  $\frac{{\rm{i}}-1}{2}$ is congruent with $i$ under the action of $\Gamma(1)$ : $j\left(\frac{{\rm{i}}-1}{2}\right)=j({\rm{i}})=1728$ which coincides with value of $1728.f(z)=64\frac{(z+3)^3}{(z-1)^2}$ at $z=0$. So  $j\left(\frac{\sqrt{2}{{\rm{i}}}}{2}\right)=j(\sqrt{2}{\rm{i}})$ must be $1728.f(2)=20^3$. Finally, when $N=3$, under the action of $\Gamma_0(3)$, $[z]=[\frac{-1}{3z}]$ for $z=\frac{\sqrt{3}{\rm{i}}}{3},z=\rho$. The $j$-function and $27\frac{z(z+8)^3}{(z-1)^3}$ vanish at $\rho$ and $z=0$ respectively. Therefore $1728.f(2)=2\times 30^3$ equals $j\left(\frac{\sqrt{3}{\rm{i}}}{3}\right)=j(\sqrt{3}{\rm{i}})$. In conclusion:
\begin{cor}\label{j-values}
\begin{equation*}
\begin{split}
&j(\sqrt{3}{\rm{i}})=16\times 15^3\quad j(2{\rm{i}})=66^3\quad j(\sqrt{2}{\rm{i}})=20^3\quad j\left(\frac{1+3\sqrt{3}{\rm{i}}}{2}\right)=-3\times 160^3\\
& j(3{\rm{i}})=\left(18+6\sqrt{3}\right)(82+54\sqrt{3})^3\quad j\left(\frac{1+3{\rm{i}}}{2}\right)=\left(18-6\sqrt{3}\right)(82-54\sqrt{3})^3
\end{split}
\end{equation*}
\end{cor}

\vskip 0.4 cm
{\bf{Acknowledgement.}} The authors would like to thank  Professors A. Rajaei and R. Takloo-Bighash for  reading the first draft of this manuscript and their valuable suggestions.

%%%%%%%%%%%%%%%%%%%%%%%%%%%%%%%%%%%%%%%%%%%%%%%%%%%%%%%%%%%%%%%%%%%%%%%%%%%%%%%%%%%%%%%%%%%%%%%%%%%%%%%%

\end{document}